\documentclass[journal,twocolumn]{IEEEtran}
\IEEEoverridecommandlockouts          % This command is only needed if you want to use the \thanks command
\usepackage[usenames,dvipsnames]{xcolor}
\usepackage{soul}

\usepackage{amsmath}
\usepackage{amssymb}
\usepackage{mathtools}
\usepackage{graphicx}
\usepackage[noadjust]{cite}

\def\str#1{{\mathcal #1}}

\usepackage{color}
\usepackage[usenames,dvipsnames]{xcolor}

\setstcolor{red}

\newcommand{\T}{\textcolor{red}{T}}
\newcommand{\R}{\textcolor{VioletRed}{R}}
\renewcommand{\P}{\textcolor{Blue}{P}}
\renewcommand{\S}{\textcolor{Cerulean}{S}}

\newtheorem{theorem}{Theorem}
\newtheorem{lemma}{Lemma}
\newtheorem{proposition}{Proposition}

\newtheorem{corollary}{Corollary}

\def\str#1{{\mathcal #1}}

\newcommand{\bd}{\mathrm{bd}}

\renewcommand{\int}{\mathrm{int}}
\renewcommand{\L}{{\str L}}

\title{\LARGE \bf Global Convergence for Replicator Dynamics of Repeated Snowdrift Games}
\author{Pouria Ramazi and Ming Cao%
\thanks{The work was supported in part by the European Research Council (ERC-StG-307207) and the Netherlands Organization for Scientific Research (NWO-vidi-14134).}
\thanks{P. Ramazi is with 
	Statistical and Mathematical Sciences Department,
	University of Alberta,
	Canada
and M. Cao is with
        ENTEG, 
	Faculty of Science and Engineering,
        University of Groningen, 
	The Netherlands,
          {\tt\small p.ramazi@gmail.com,\,m.cao@rug.nl }
    }
}
\begin{document}

\maketitle
\thispagestyle{empty}
\pagestyle{empty}
\begin{abstract}
	To understand the emergence and sustainment of cooperative behavior in interacting collectives, we perform global convergence analysis for replicator dynamics of a large, well-mixed population of individuals playing a repeated snowdrift game with four typical  strategies, which are always cooperate (ALLC), tit-for-tat (TFT), suspicious tit-for-tat (STFT) and always defect (ALLD). 
	The dynamical model is a three-dimensional ODE system that is parameterized by the payoffs of the base game. 
	Instead of routine searches for evolutionarily stable strategies and sets, we expand our analysis to determining the asymptotic behavior of solution trajectories starting from any initial state, and in particular show that for the full range of payoffs, every trajectory of the system converges to an equilibrium point. 
	What enables us to achieve such comprehensive results is studying the dynamics of \b{two ratios} of the state variables, each of which either monotonically increases or decreases in the half-spaces separated by their corresponding planes.  	
	The convergence results highlight three findings that are of particular importance for understanding the cooperation mechanisms among self-interested agents playing repeated snowdrift games. 
	First, the inclusion of TFT- and STFT-players, the two types of conditional strategy players in the game, increases the share of cooperators of the overall population compared to the situation when the population consists of only ALLC- and ALLD-players.
	This confirms findings in biology and sociology that reciprocity may promote cooperation in social collective actions, such as reducing traffic jams and division of labors, where each individual may gain more to play the opposite of what her opponent chooses.
	Second, surprisingly enough, regardless of the payoffs, there always exists a set of initial conditions under which ALLC players do not vanish in the long run, which does not hold for all the other three types of players. 
	So an ALLC-player, although perceived as the one that can be easily taken advantage of in snowdrift games, has certain endurance in the long run. 
	Third, the parametric framework makes it possible to actually control the final population shares, a challenging topic in population dynamics, by tuning the payoffs of the base game. 
\end{abstract}

\section{Introduction}
	Game theory provides a framework for studying various control problems such as robust control, distributed control and optimization for traffic systems, communication networks and multi-agent systems in general; in this context, the different types of games that have been modeled and analyzed in the literature include potential games \cite{marden2012price, cortes2015self, li2013designing, marden2014generalized, li2014decoupling},
stochastic games \cite{altman2010markov,chang2003two, wikecek2011stochastic},
constrained games \cite{altman2009constrained},
repeated games \cite{marden2009joint, shamma2005dynamic},
matrix games \cite{bopardikar2013randomized},
networked games \cite{guo2013algebraic},
and others \cite{frihauf2012nash, stankovic2012distributed, vamvoudakis2014detection, marden2012state, mylvaganam2015constructive, gharesifard2012evolution, ramazi2017asynchronous}.
More recently, evolutionary game theory has gained more attention since it is a powerful tool in understanding the evolution of cooperation among selfish individuals as reported by biologists, sociologists, economists, etc \cite{van2015focus,nowak2006evolutionary,sandholm2010population,weibull1997,ramazi2015feeling, ramazi2016networks}. 
Researchers have found that network topology \cite{rand2014static}, phenotypic interactions \cite{jansen2006altruism, ramazi2016evolutionary}, punishment \cite{van2012social}, population heterogeneity \cite{ramazi2015analysis}, as well as other components in game setups can all affect the success of cooperators in face of defectors. 
One stimulating mechanism for the evolution of cooperation that is generally believed to promote cooperation, especially in human societies \cite{van2012direct}, is \emph{direct reciprocity} \cite{trivers1971evolution}.
This mechanism is captured by \emph{repeated games} where individuals play a base game repeatedly and can base their action in each round of the game on that of the opponent in the previous round, resulting in \emph{reactive strategies}.

Perhaps the most typical reactive strategy is the simple yet successful \emph{tit-for-tat ($TFT$)} strategy where the player starts with cooperation and cooperates if the opponent cooperated and defects if the opponent defected in the last round.
A more defective version of the strategy is the \emph{suspicious tit-for-tat ($STFT$)} strategy which is the same as $TFT$ except that the player starts with defection.
In addition to these conditional strategies, there are two unconditional ones which are the two extreme strategies in repeated 2-strategy games: \emph{always-cooperate ($ALLC$)} and \emph{always-defect ($ALLD$)}. 
While much research has been carried out to investigate the performance of different reactive strategies under the prisoner's dilemma game, the cornerstone of game theory, \cite{nowak1990evolution, ioannou2014asymptotic, hilbe2013evolution, grujic2012three, lorberbaum1994no, imhof2005evolutionary}, less has been devoted to the anti-coordination snowdrift game \cite{qi2015experiments, kummerli2007human,chen} despite the fact that the snowdrift game  captures many behavioral patterns that cannot be well-modeled by the prisoner's dilemma game \cite{hauert2004spatial}.
Moreover, the existing results on the snowdrift game are mainly experimental or simulation based. 
For example, in \cite{kummerli2007human}, based on human experiments the authors postulate that iterated snowdrift games can explain high levels of cooperation among non-relative humans. 
However, few mathematical statements have been constructed to support such claims \cite{doebeli2004evolutionary, hauert2006synergy, madeo2014game, qin2017neighborhood}.

The performance of different reactive strategies also remains an open problem. Usually the strategies are compared using 2-strategy games, e.g., the two famous competitions conducted by Axelord \cite{axelrod1980effective, axelrod1980more} where strikingly, the simple $TFT$ was placed first in both (note that although $TFT$ is known to be successful mostly in the repeated prisoner's dilemma, it has also been reported to be successful in the repeated snowdrift game \cite{kummerli2007human, dubois2003forager}). 
The situation would be different if more than two strategies could be played in the game. 
Then the best strategy can be decided by natural selection, which is captured by evolutionary dynamics such as the well-known \emph{replicator dynamics} \cite{ben2014delayed, borkar2003dynamic, brunetti2015state, drighes2014stability}. 
Due to nonlinearity, the replicator dynamics, however, may exhibit quite complex behaviors, as the dynamically-equivalent \emph{Lotka-Volterra Equations} do \cite{hofbauer1998evolutionary}.  
Indeed, except for a few cases \cite{diekmann2009cyclic}  \cite{zeeman2003local}, the analysis is restricted to only those modeled by planar dynamical systems \cite{Bomze}. 
This makes the performance investigation of more than three reactive strategies generally challenging under the replicator dynamics. 
However, the assumption of having just a small number of available strategies may seem not to be realistic or representative for many natural phenomena, particularly those involving a wide range of mutations taking place. 
A research line has consequently been established to study evolutionary outcomes of repeated games with a large or possibly infinite number of reactive strategies by limiting the analysis to finding \emph{evolutionarily stable strategies} and \emph{sets}, which are known to be asymptotically stable under many evolutionary dynamics such as the replicator dynamics \cite{weibull1997}. 
For example, the repeated prisoner's dilemma is shown to have no pure strategies that are evolutionarily stable or that can form an evolutionarily stable set \cite{selten1984gaps, garcia2016and}.
Although revealing (non)existence of stable sets under the evolutionary dynamics, these works neglect other possible long-run behaviors, such as a saddle point as the simplest example.
Thus, a considerable portion of equilibrium states that can be favored by natural selection remains concealed. 
Moreover, having many available reactive strategies is not always a reasonable assumption, especially when complex strategies are costly or uncommon \cite{samuelson2003evolutionary}.
So there is a need for exhaustive asymptotic analysis of evolutionary dynamics with typical and simple reactive strategies.  
The convergence of large populations playing evolutionary games is of general interest and has applications in control theory  \cite{theodorakopoulos2012selfish, madeo2014game, obando2013building, wiecek2010stochastic}.

We address both of the above issues in this paper. 
While considering the snowdrift game as the base game, we study the evolution of a large population of individuals playing the four just mentioned strategies, $ALLC$, $TFT$, $STFT$ and $ALLD$, under the replicator dynamics.
We consider a completely parameterized payoff matrix with an arbitrary number of repetitions of the base game and reveal all asymptotic outcomes of the resulting 3-dimensional dynamics.
What enables us to expand our analysis beyond the routine search for evolutionarily stable sets is studying the dynamics of two ratios of the state variables. 
By dividing the simplex into four sections, in each of which each ratio either monotonically increases or decreases, we show that every trajectory of the system converges to an equilibrium point, excluding the possibility of limit cycles or chaotic behaviors. 
This approach can be applied to general replicator dynamics with more than three strategies where one or more ratios of the state variables monotonically increase or decrease in some part of the simplex. 
Our analyses shed light on the social dilemma in the snowdrift game, that is why selfish individuals cooperate while they earn more if they defect against their cooperative opponents. This is done by showing that first of all, even in the presence of the very defective strategy $ALLD$, for some range of payoffs and initial population portions, the population evolves to the state where all mutually cooperate. In other words, natural selection disfavors individuals playing $ALLD$ and instead chooses those playing more cooperative strategies such as $TFT$ and even $ALLC$.
Secondly, the convergence results postulate that among the four types of players, $ALLC$s are surprisingly the best in terms of survival and appearance in the long run, explaining why selfish individuals may repeatedly cooperate in a snowdrift social context. 
As a second contribution, due to the parametric framework we provide, our convergence analysis can be used to actually control the final state of the replicator dynamics. 
By tuning the parameters, one can control the final population portions of individuals playing the reactive strategies. 
This is possible when a central agency has control over the payoff matrix, e.g., tax regulations made by the government.
Moreover, for populations initially having four co-existing types of players, by comparing those final states in which one or two types of players die out to those with all four, it becomes clear how adding a third or fourth strategy can change the final population state. These results lead to addressing the crucial question of \emph{how to control portions of different types of individual in a decision-making population?}, which finds fascinating applications in repeated snowdrift games, ranging from trading commodities to division of labor.

The rest of the paper is organized as follows.
In Section \ref{sec-pf}, we describe the replicator dynamics for repeated snow drift games with the above four reactive strategies. 
In Section \ref{secGlobal}, we provide the global convergence results and discuss their implications on the success of the strategies. 
We end with the concluding remarks in Section \ref{secConclusion}.

\section{Problem formulation}		\label{sec-pf}
We consider an infinitely large, well-mixed population of
individuals that are playing repeated games over time. Each game has
two players with two pure strategies: one is to cooperate, denoted
by $C$, and the other to defect, denoted by $D$, and the payoffs of
the game, described by the following payoff matrix, are symmetric to
both players
\begin{equation}                                                    \label{baseGame}
    \bordermatrix{
        ~   &   C   &   D   \cr
        C   &   \R  &   \S  \cr
        D   &   \T  &   \P  \cr
    },
\end{equation}
where $\R$, $\S$, $\T$ and $\P$ are real numbers and sometimes in the
literature are called the reward, sucker's payoff, temptation and
punishment respectively. 
We call this two-player,
symmetric game, the \emph{base game} and denote it by $G$. When
the payoffs of the game satisfy
\begin{equation}                                                    \label{snowdriftInequality}
    \T>\R>\S>\P,
\end{equation}
the game is called a \emph{snowdrift game} (also known as the \emph{hawk-dove} or the \emph{chicken game}).
The game has two Nash equilibria in pure strategies, both of which correspond to the situation when the two players play
opposite strategies, and for this reason such a game is also called
an \emph{anti-coordination} game, often used to study how players
may contribute to the accomplishment of a common task. 
In this study, we are particularly interested in the case in which individuals play the game repeatedly over time and adjust their
strategies according to what their opponents have played in the
past. 
Formally, a \emph{repeated game, denoted by $G^{m}$}, $m\geq
2$, with \emph{reactive strategies} is constructed from the base
game $G$ by repeating it for $m$ rounds, and limiting a player's
choice of strategies in the current round to be based on the
opponent's choice in the previous round. 
In fact, a reactive strategy $s$ can always be represented by the triple $(p,q,r)$,
where $p$ is the probability of cooperating in the first round, and $q$ (respectively $r$) is the probability of cooperating if the opponent has cooperated (respectively defected) in the previous round. We consider the following  strategies:
\begin{itemize}
    \item \emph{always-cooperate (ALLC)}, $(1,1,1)$: always cooperates;
    \item \emph{tit-for-tat (TFT)},  $(1,1,0)$: cooperates in the first round, and then chooses what the opponent did in the previous round;
    \item \emph{suspicious-tit-for-tat (STFT)},  $(0,1,0)$: defects in the first round, and then chooses what the opponent did in the previous round;
    \item \emph{always-defect (ALLD)}, $(0,0,0)$: always defects.
\end{itemize}
When two players play the repeated game $G^m$, the payoffs for the reactive strategies
can be calculated every $m$ rounds, leading to the payoff matrix $A:=[a_{ij}]$ defined by
%following $4\times4$ payoff matrix
\begin{gather*}
	A  =  	
\scalebox{.66}{
		\bordermatrix{
~	&ALLC 			& TFT 										& STFT										&ALLD		\cr
ALLC&m\R  			&m\R										&\S+(m-1)\R									&m\S		\cr
TFT 	&m\R 			&m\R										&\lceil\frac{m}{2}\rceil\S+\lfloor\frac{m}{2}\rfloor\T 	&\S+(m-1)\P\cr
STFT	&\T +(m-1)\R 	&\lceil\frac{m}{2}\rceil\T+\lfloor\frac{m}{2}\rfloor\S 	&m\P										&m\P		\cr
ALLD&m\T			&\T + (m-1)\P									&m\P										&m\P		\cr
	}
}.
\end{gather*}
To illustrate how the matrix $A$ is obtained, we take the match between $TFT$ and $ALLD$ as an example.
In round one, the $TFT$ player cooperates and the $ALLD$ player defects, so their payoffs according to (\ref{baseGame}) are $\S$ and $\T$ respectively. 
From round two, both $TFT$ and $ALLD$ players defect and hence receive $\P$. 
So over time the payoffs for the $TFT$ player are $\S, \P,  \P, \ldots$ while those for the $ALLD$ player are $\T, \P, \P, \ldots$. 
Summing up the payoffs over the $m$ rounds, one obtains the entries of $a_{23}$
and $a_{32}$ in $A$. 
Hence, the repeated game $G^m$ can be taken as
a normal, symmetric two-player game with the payoff matrix $A$ and
with the pure-strategy set $\{ ALLC, TFT, STFT, ALLD\}$.

Restricting the base game to be played $m$ rounds with the same opponent is an assumption that holds in many natural systems and real-life scenarios. 
Birds in the same flock migrating to winter quarters interact with each other during periods of their migration; students in the same project group collaborate with each other during the semester; and tenants in the same apartment meet each other during their rental period.
Such interactions take place repeatedly with the same individuals for a certain amount of time.

Having clarified how a pair of individuals play games with each other, we now
describe the evolutionary dynamics of the whole population. 
Towards this end, we introduce \emph{replicator dynamics}, which is a standard model from evolutionary game theory \cite{weibull1997,sandholm2010population}.

Let $0\leq x_i(t) \leq 1, i=1,2,3$ and $4$, denote the population shares at time $t$ of those
individuals playing the pure strategies $ALLC$, $TFT$, $STFT$ and $ALLD$
respectively. 
Since the four types of players constitute the whole population, it follows that for all $t$, $\sum_{i=1}^4 x_i= 1$.
Define the population vector $$x := \begin{bmatrix} x_1 &  x_{2} &  x_3 &  x_{4} \end{bmatrix} ^\top .$$
Then $x \in \Delta$ where $\Delta$ is the 4-dimensional simplex defined by
\begin{equation}	\label{delta}
	\Delta: = \left\{ z \,|\, z\in\mathbb{R}^4, z_i \geq 0, i=1,\ldots,4, \sum_{i=1}^4 z_i = 1\right\}.
\end{equation}
We use the unit vectors at the vertices of the simplex
\begin{equation*}
    p^1  = \begin{bmatrix} 1 \\ 0 \\ 0 \\ 0 \end{bmatrix},
    p^2  = \begin{bmatrix} 0 \\ 1 \\ 0 \\ 0 \end{bmatrix},
    p^3  = \begin{bmatrix} 0 \\ 0 \\ 1 \\ 0 \end{bmatrix},    
    p^4  = \begin{bmatrix} 0 \\ 0 \\ 0 \\ 1 \end{bmatrix}
\end{equation*}
to represent the population vectors corresponding to all $ALLC$ players, all $TFT$ players, all $STFT$ players and all $ALLD$ players respectively.
Then the evolutions of  
$x_i$, $i=1,\ldots,4$, are described by the replicator dynamics \cite{weibull1997,sandholm2010population}
\begin{equation}
    \dot{x}_i = [u(p^i,x) - u(x,x)]x_i,                         \label{RD}
\end{equation}
where $u(\cdot,\cdot)$ is the \emph{utility function} defined by
\[
    u(x,y)=x^\top  A y \textrm{\ for\ } x,y \in \Delta
\]
determining the fitness of a player.
In essence, \eqref{RD} indicates that in an evolutionary process, the reproduction rate of the strategy-$i$ players is proportional to the difference between the fitness of strategy-$i$ players $u(p^i,x)$ and the average population fitness $u(x,x)$ as a consequence of the fact that the more payoff an individual acquires when playing against its opponents, compared to the average payoff of the whole population, the more new offspring proportionally it produces. 
Since $u(p_i,x)$ is the expected payoff of an $i$-playing individual against a random other individual in the population, the dynamics can be shown to be interpretable as follows \cite{weibull1997}.  
Over a continuous course of time, an individual in the population (say a \emph{TFT} player) randomly meets another (say an \emph{ALLD} player), plays the base game with her opponent for $m$ rounds, earns an accumulated payoff according to the payoff matrix $A$ (that is $S+(m-1)P$), and reproduces offspring playing her same strategy (these are \emph{TFT} players) with a rate equal to her payoff. 
Indeed, there are two time scales of fast and slow dynamics; the time it takes for two players to play the repeated game goes much faster and in fact neglectable compared to the reproduction time.
The dynamics can also be seen as the mean dynamic approximation of the following process that takes place over a discrete sequence of time \cite{sandholm2010population}. 
At each time step, 
\emph{i)} every individual plays the base game for $m$ rounds with every other individual in the population and earns the payoff of the average, and 
\emph{ii)} a random individual updates her reactive strategy according to the \emph{pairwise proportional imitation update rule}, that is, she randomly chooses another individual, say $j$, and if her payoff is less than that of individual $j$, imitates his strategy with a probability proportional to the payoff difference, and otherwise, sticks to her own current strategy.

Since $u(\cdot,\cdot)$ is continuously differentiable in $\mathbb{R}^4\times \mathbb{R}^4$, \eqref{RD} has a unique solution for any $x(0) \in\Delta$ \cite[Theorem  7.1.1]{wiggins2003introduction}.
It is easy to check that the solution indeed satisfies the constraints $0\leq x_i(t) \leq 1$, $i=1,\ldots,4$ for all $t$.
Moreover, it can be verified that for any $t$, if $x(t)\in\Delta$, it holds that $\sum_{i=1}^4 \dot{x}_i (t) = 0$. 
Hence, $\sum_{i=1}^4 x_i= 1$ is in force for all $t$ given $x(0)\in\Delta$. 
%proof: let z(t) = \sum_{i=1}^3 x_i(t). For any $t_1$, if $z(t_1)=1$, $\dot{z}(t_1) = 0$. Hence, $z(t) = 1$ for all $t$. Also contradiction can be used.  
Therefore, $\Delta$ is invariant under \eqref{RD} and hence the dynamical system \eqref{RD} is well defined on $\Delta$.
% see proposition 6.1 in Weibull.

We perform global convergence analysis of the replicator dynamics \eqref{RD}. More specifically, for any given initial condition $x(0)\in\Delta$, we aim to determine the limit state of $x(t)$ for \eqref{RD}.

%%%
\section{Global convergence result}	\label{secGlobal}
The main results of this paper are presented in this section. 
First we find the equilibrium points of the system.
Then for the convergence results, we divide the analysis into several parts using the notion of \emph{face} defined as follows.
A {face} of the simplex is the convex hull of a non-empty subset $\str H$ of $\{p^1,p^2,p^3,p^4\}$, and is denoted by $\Delta(\str H)$. 
For simplicity, we remove the braces when $\str H$ is represented by its members.
For example, the face $\Delta(p^1,p^3,p^4)$ is the convex hull of $\str H = \{p^1,p^3,p^4\}$.
When $\str H$ is proper, $\Delta(\str H)$ is called a \emph{boundary face}.
Following convention, the \emph{boundary} of a set $\str S$, denoted by $\mathrm{bd}(\str S)$, is the set of points $p$ such that every neighborhood of $p$ includes at least one point in $\str S$ and one point out of $\str S$, and the \emph{interior} of  $\str S$, denoted by $\mathrm{int}(\str S)$, is the greatest open subset of $\str S$. 
The following result enables us to analyze the evolution of a trajectory starting from $\bd(\Delta)$ separately from that starting from $\int(\Delta)$.

\begin{lemma}	\label{faceIsInvariant}
	Each face of $\Delta$ is invariant under the replicator dynamics \eqref{RD}.
\end{lemma}
\begin{IEEEproof}
	$\Delta$ was already shown to be invariant in Section \ref{sec-pf}.
	So it remains to prove the lemma for the boundary faces. 
	This can be done based on the observation that if for some $i=1,\ldots,4$, $x_i(0) = 0$,
	then $x_i(t) = 0$ for all $t\in\mathbb{R}$.	
\end{IEEEproof}

We start with analyzing the boundary of the simplex.
However, the boundary of the simplex itself consists of the four planar faces $\Delta(p^1,p^2,p^3)$, $\Delta(p^1,p^2,p^4)$, $\Delta(p^1,p^3,p^4)$ and $\Delta(p^2,p^3,p^4)$.
Because of Lemma \ref{faceIsInvariant}, we can also analyze the dynamics \eqref{RD} on each of these faces separately. 
Yet again, the boundary of each of these planar faces consists of three one-dimensional faces known as the \emph{edges} of the simplex. 
For example, the boundary of the face $\Delta(p^1,p^2,p^3)$ consists of the edges $\Delta(p^1,p^2)$, $\Delta(p^1,p^3)$ and $\Delta(p^2,p^3)$.
On the other hand, each of the edges are also invariant in view of Lemma \ref{faceIsInvariant}.
Therefore, we study separately trajectories starting from an edge and those  starting from the interior of a planar face.
Then we proceed to the interior of the simplex.

To simplify the analysis, we carry out on the matrix $A$ some operations that preserve the dynamics (\ref{RD}). 
Subtracting $m\R$ from the entries of the first and second columns, and $m\S$ from the entries of the third and fourth columns of $A$, we acquire the following matrix
	\begin{gather}	
		A' :=	[a'_{ij}] = 	\nonumber\\
		\scalebox{.72}{$		
		\begin{bmatrix}
			0	 	&0																&\S+(m-1)\R-m\P											&m(\S-\P)	\\
 			0	 	&0																&\lceil\frac{m}{2}\rceil\S+\lfloor\frac{m}{2}\rfloor\T  - m\P		 	&\S-\P		\\
			\T -\R	&\lceil\frac{m}{2}\rceil\T+\lfloor\frac{m}{2}\rfloor\S - m\R					&0														&0			\\
			m(\T-\R)	&\T + (m-1)\P	-m\R												&0														&0			\\
		\end{bmatrix}
		$}.		\label{A'}
	\end{gather}
In view of Lemma \ref{lem-RDInvariantUnderAddition} in Appendix \ref{app:A}, the dynamics \eqref{RD} are unchanged with $A'$ in place of $A$.
Since $A'$ is more structured with zero block matrices, in what follows we focus  on $A'$ instead of $A$. 
\subsection{Equilibrium points} 
To determine the equilibria of the system, we first look for those on the boundary of the simplex, and then for those in the interior.

\subsubsection{boundary equilibrium points}
Let $\Delta^o$ and $\Delta^{oo}$ denote the set of equilibrium points of the replicator dynamics \eqref{RD} that belong to $\Delta$ and $\mathrm{bd}(\Delta)$, respectively. 
Depending on the payoffs, $\Delta^{oo}$ will be a combination of the unit vectors $p^1,p^2,p^3, p^4$, the vectors
\begin{gather*}
	x^{14}  =	
	\begin{bmatrix} 	
		\frac{\S-\P}{\S-\P+\T-\R}	\\	0	\\	0	\\	\frac{\T-\R}{\S-\P+\T-\R}	
	\end{bmatrix},
	x^{23} = 
	\begin{bmatrix} 	
		0 \\ \frac{ \lceil\frac{m}{2}\rceil\S+\lfloor\frac{m}{2}\rfloor\T  - m\P	}{m(\T+\S-\P-\R)}	\\	\frac{\lceil\frac{m}{2}\rceil\T+\lfloor\frac{m}{2}\rfloor\S - m\R	}{m(\T+\S-\P-\R)}	\\	0	
	\end{bmatrix},		\\
	x^{13} = 
	\begin{bmatrix}	
		\frac{\S+(m-1)\R-m\P}{\T+\S+(m-2)\R -m\P}	\\	0\\	\frac{\T-\R}{\T+\S+(m-2)\R-m\P}	\\	0
	\end{bmatrix},
	x^{24} = 
	\begin{bmatrix}
		0	\\	\frac{\S-\P}{\T+\S+(m-2)\P-m\R}	\\	0	\\	\frac{\T+(m-1)\P-m\R}{\T+\S+(m-2)\P-m\R}
	\end{bmatrix},
\end{gather*}
and the sets
\begin{align*}	
	\str X^{12} &= \{ \alpha p^1 + (1-\alpha)p^2 : \alpha \in [0,1] \},		\\
	\str X^{34} &= \{ \alpha p^3 + (1-\alpha)p^4 : \alpha \in [0,1] \},		\\
	\str X^{123} &= \{ x\in \mathrm{int}(\Delta(p^1,p^2,p^3))\, |\, a'_{31}x_1+a'_{32}x_2-a'_{13}x_3 = 0\}
\end{align*}		
where $a'_{ij}$'s are the entries of $A'$ defined in \eqref{A'}.
Here, the superscript $ij$ in $x^{ij}$ (resp. $\str X^{ij}$) simply means that $x^{ij}$ (resp. $\str X^{ij}$) belongs to the edge $\Delta(p^i,p^j)$.
The following proposition determines $\Delta^{oo}$.
\begin{proposition}	\label{lem-equilibria}
	Assume \eqref{snowdriftInequality} holds. 
	It follows that
	\begin{enumerate}
		\item if $\S<\R< \frac{\T+(m-1)\P}{m}$, then 
		 \begin{equation*}
		 	\Delta^{oo} =  \str X^{12} \cup \{x^{13}, x^{14}, x^{23},x^{24}\} \cup  \str X^{34};
		 \end{equation*}	
		 
		 \item if $\frac{\T+(m-1)\P}{m} \leq \R< \frac{\T+\S}{2}$, or if $m=2n+1,n\geq1$ and $\frac{\T+\S}{2} < \R < \frac{(n+1)\T+n\S}{2n+1}$, then 
		 \begin{equation*}
		 	\Delta^{oo} = \str X^{12} \cup \{x^{13}, x^{14}, x^{23}\} \cup  \str X^{34};
		 \end{equation*}		
		 
		 \item if $m=2n+1,n\geq1$ and $\R=\frac{\T+\S}{2}$, then 
		 \begin{equation*}
		 	\Delta^{oo} =\str X^{12} \cup \{x^{13}, x^{14}, x^{23}\} \cup  \str X^{34}\cup \str X^{123};
		 \end{equation*}	
		 
		  \item if $m=2n,n\geq1$ and $\R = \frac{n\T+(n-1)\S}{2n-1}$, then 
		 \begin{equation*}
		 	\Delta^{oo} =\str X^{12} \cup \{x^{13}, x^{14}\} \cup  \str X^{34}\cup \str X^{123};
		 \end{equation*}	
		
		\item if $\max\left\{ \frac{\lceil\frac{m-2}{2}\rceil\S+\lfloor\frac{m}{2}\rfloor\T}{m-1} ,\frac{\lceil\frac{m}{2}\rceil\T+\lfloor\frac{m}{2}\rfloor\S}{m}    \right\}< \R < \T$, or if $m=2n,n\geq1$ and $\frac{\T+\S}{2} \leq \R< \frac{n\T+(n-1)\S}{2n-1}$, then 
		 \begin{equation*}
		 	\Delta^{oo} = \str X^{12} \cup \{ x^{13}, x^{14}\} \cup  \str X^{34}.
		 \end{equation*}		 		 
	\end{enumerate}
\end{proposition}

For the proof,  we need to take a closer look at the payoff matrix $A'$.
The order in the magnitudes of the entries in each column of $A'$, clarified in the following lemma, proves useful both in the determination of the equilibria and the asymptotic behavior of the replicator dynamics \eqref{RD}.
\begin{lemma}		\label{lem-signStructure}
	Assume \eqref{snowdriftInequality} holds. 
	Consider the payoff matrix $A'$ and denote the maximum positive, positive, negative and minimum negative entries of each column by `$++$', `$+$' , `$-$'  and `$--$', respectively. 
	Then $A'$ has the following sign structure
	\begin{enumerate}
		\item \scalebox{.85}{$
		\begin{bmatrix}
			0	 	&0		&+		&++		\\
 			0	 	&0		&++	 	&+		\\
			+		&++		&0		&0		\\
			++		&+		&0		&0		\\
		\end{bmatrix}$}
		when $\S<\R < \tfrac{\T+(m-1)\P}{m};$		

		\item \scalebox{.85}{$
		\begin{bmatrix}
			0	 	&0		&+		&++		\\
 			0	 	&0		&++	 	&+		\\
			+		&++		&0		&0		\\
			++		&0,-		&0		&0		\\
		\end{bmatrix}$}
		when $\tfrac{\T+(m-1)\P}{m} \leq \R <  \tfrac{\T+\S}{2};$

		\item \scalebox{.85}{$
		\begin{bmatrix}
			0	 	&0		&++		&++		\\
 			0	 	&0		&++	 	&+		\\
			+		&++		&0		&0		\\
			++		&-		&0		&0		\\
		\end{bmatrix}$}
		when $m=2n+1,n\geq 1, \text{ and } \R = \tfrac{\T+\S}{2}; $	
			
		\item \scalebox{.85}{$
		\begin{bmatrix}
			0	 	&0		&++		&++		\\
 			0	 	&0		&+	 	&+		\\
			+		&++,0	&0		&0		\\
			++		&-		&0		&0		\\
		\end{bmatrix}$}		
		when $m=2n+1,n\geq 1, \text{ and } \tfrac{\T+\S}{2} < \R \leq \tfrac{(n+1)\T+n\S}{2n+1};$

		\item \scalebox{.85}{$
		\begin{bmatrix}
			0	 	&0		&+		&++		\\
 			0	 	&0		&++	 	&+		\\
			+		&0,-		&0		&0		\\
			++		&--		&0		&0		\\
		\end{bmatrix}$}
		when $m=2n,n\geq 1, \text{ and } \tfrac{\T+\S}{2} \leq \R < \tfrac{n\T+(n-1)\S}{2n-1};$ 		
	
		\item \scalebox{.85}{$
		\begin{bmatrix}
			0	 	&0		&++		&++		\\
 			0	 	&0		&++	 	&+		\\
			+		&-		&0		&0		\\
			++		&--		&0		&0		\\
		\end{bmatrix}$}
		when $m=2n,n\geq 1,  \text{ and } \R = \tfrac{n\T+(n-1)\S}{2n-1}; $
	
		\item and \scalebox{.85}{$ 
		\begin{bmatrix}
			0	 	&0		&++		&++		\\
 			0	 	&0		&+	 	&+		\\
			+		&-		&0		&0		\\
			++		&--		&0		&0		\\
		\end{bmatrix}$}		
		when $\max\left\{ \tfrac{\lceil\frac{m-2}{2}\rceil\S+\lfloor\frac{m}{2}\rfloor\T}{m-1} ,\tfrac{\lceil\frac{m}{2}\rceil\T+\lfloor\frac{m}{2}\rfloor\S}{m}    \right\}<\R<\T.$
	\end{enumerate}
	Here, when an entry takes both $0$ and one other sign (separated by a comma), $0$ takes place if the equality sign of the $\R$ condition holds, and otherwise the other sign is valid. 
\end{lemma}
\begin{IEEEproof}
	The sign of the elements of $A'$ are determined by \eqref{snowdriftInequality}.
	First note that $\T>\R$ implies $a'_{31}>0$. 
	On the other hand, since $m\geq2$, we have that $a'_{41}>a'_{31}>0$. 
	Hence, due to the fact that the third and fourth entries of the last column of $A'$ are zero, $a'_{41}$ and $a'_{31}$ are denoted by `$++$' and `$+$', respectively.
	Similarly $\S>\P$ implies $a'_{14} > a'_{24}>0$ and hence $a'_{14}$ and $a'_{24}$ are denoted by `$++$' and `$+$', respectively.
	Since $\T,\S>\P$ implies
	\begin{equation*}
		\lceil\frac{m}{2}\rceil\S+\lfloor\frac{m}{2}\rfloor\T  
		> \lceil\frac{m}{2}\rceil\P+\lfloor\frac{m}{2}\rfloor\P  	 
		\Rightarrow
		 \lceil\frac{m}{2}\rceil\S+\lfloor\frac{m}{2}\rfloor\T 	
		> m\P ,								
	\end{equation*}
	it follows that $a'_{23}>0$.
	Additionally, 
	$$
		\R > \P, \ \S > \P 
		\Rightarrow 
		(m-1)\R + \S >m\P,   
	$$
	which implies $a'_{13}>0$. 
	Similarly $\T,\S > \P$ yields
	\begin{gather*}
	\scalebox{.97}{$	
		\T + \lceil\frac{m-2}{2}\rceil\T+\lfloor\frac{m-2}{2}\rfloor\S + \S 
		> 
		\T + \lceil\frac{m-2}{2}\rceil\P+\lfloor\frac{m-2}{2}\rfloor\P + \P	$}	\\
	\scalebox{0.97}{$		
		\Rightarrow
		 \lceil\frac{m}{2}\rceil\T+\lfloor\frac{m}{2}\rfloor\S 
		 >
		 \T + (m-2)\P + \P 
		 = \T + (m-1)\P
	$}. 
	\end{gather*}
	Hence, $a'_{32}>a'_{42}$.
	It remains to determine the signs of $a'_{42}$ and $a'_{32}$ and also the ordering of $a'_{13}$ and $a'_{23}$.
	Since $m\geq 2$, division by $m-1$ is valid, and hence the following hold
	\begin{gather}
		a'_{42} > 0 
		\iff 
		\R < \frac{\T+(m-1)\P}{m},		\label{lem-signStructure-1}\\
		a'_{32} > 0
		\iff 
		\R < \frac{\lceil\frac{m}{2}\rceil\T+\lfloor\frac{m}{2}\rfloor\S}{m}, 	\label{lem-signStructure-2}\\
		a'_{23} > a'_{13} 
		\iff
		\R < \frac{\lceil\frac{m-2}{2}\rceil\S+\lfloor\frac{m}{2}\rfloor\T}{m-1}	.	\label{lem-signStructure-3}	
	\end{gather}
	
	The average of $\T, \underbrace{\P, \ldots, \P}_{m-1}$ is less than both the average of $\T, \underbrace{\T,\ldots,\T}_{\lceil\frac{m}{2}\rceil-1}, \underbrace{\S, \ldots, \S}_{\lfloor\frac{m}{2}\rfloor}$ and the average of $\T, \underbrace{\T,\ldots,\T}_{\lfloor\frac{m}{2}\rfloor-1}, \underbrace{\S, \ldots, \S}_{\lceil\frac{m-2}{2}\rceil}$. 
	Thus, 
	\begin{equation*}
		\frac{\T+(m-1)\P}{m}
		<
		\frac{\lceil\frac{m}{2}\rceil\T+\lfloor\frac{m}{2}\rfloor\S}{m},
		 \frac{\lceil\frac{m-2}{2}\rceil\S+\lfloor\frac{m}{2}\rfloor\T}{m-1}. 
	\end{equation*}
	Hence, when \eqref{lem-signStructure-1} holds, so do \eqref{lem-signStructure-2} and \eqref{lem-signStructure-3}.
	This proves the first case of the lemma. 
	Now we compare $ \frac{\lceil\frac{m-2}{2}\rceil\S+\lfloor\frac{m}{2}\rfloor\T}{m-1} $ and $\frac{\lceil\frac{m}{2}\rceil\T+\lfloor\frac{m}{2}\rfloor\S}{m} $.
	In general, it holds that 
	\begin{align*}
		 \dfrac{\lceil\frac{m-2}{2}\rceil\S+\lfloor\frac{m}{2}\rfloor\T}{m-1} = 
		 \left\{\begin{aligned}
		 		&\dfrac{(n-1)\S + n\T}{2n-1} 		&m=2n				\\
		 		&\dfrac{n\S + n\T}{2n} 			&m=2n +1 		
		 \end{aligned}\right. ,
	\end{align*}
	\begin{equation*}
		\dfrac{\lceil\frac{m}{2}\rceil\T+\lfloor\frac{m}{2}\rfloor\S}{m}  = 
		 \left\{\begin{aligned}
		 		&\dfrac{n\T + n\S}{2n}	 			&m=2n			\\
		 		&\dfrac{(n+1)\T + n\S}{2n+1} 		&m=2n +1 		
		 \end{aligned}\right. .
	\end{equation*}
	Due to the fact that $\T>\S$, we obtain
	\begin{equation}	\label{lem-signStructure-4}
		\begin{aligned}
		\dfrac{\lceil\frac{m}{2}\rceil\T+\lfloor\frac{m}{2}\rfloor\S}{m}
		 < \dfrac{\lceil\frac{m-2}{2}\rceil\S+\lfloor\frac{m}{2}\rfloor\T}{m-1} & \quad m=2n ,		\\
		\dfrac{\lceil\frac{m-2}{2}\rceil\S+\lfloor\frac{m}{2}\rfloor\T}{m-1}
		 < \dfrac{\lceil\frac{m}{2}\rceil\T+\lfloor\frac{m}{2}\rfloor\S}{m} & \quad m=2n+1. 		%\label{inequalityOdd}
		 \end{aligned}
	\end{equation}
	Hence, 
	\begin{equation}	\label{lem-signStructure-5}
		\min\left\{  
			\dfrac{\lceil\frac{m}{2}\rceil\T+\lfloor\frac{m}{2}\rfloor\S}{m},
		 	\dfrac{\lceil\frac{m-2}{2}\rceil\S+\lfloor\frac{m}{2}\rfloor\T}{m-1}	
		\right\}
		= \frac{\T+\S}{2}.
	\end{equation}
	The above equation results in cases 2) and 3) of the lemma. 
	The remaining cases can be verified similarly using \eqref{lem-signStructure-4} and \eqref{lem-signStructure-5}.
\end{IEEEproof}

The boundary of $\Delta$ is the union of the boundary faces $\Delta(p^1,p^2,p^3)$, $\Delta(p^1,p^2,p^4)$, $\Delta(p^1,p^3,p^4)$ and $\Delta(p^2,p^3,p^4)$. 
So in order to find the equilibria on $\bd(\Delta)$, we can investigate each face separately.  
The interior equilibria of each face is determined in the following proposition, the proof of which follows from the convergence results and methods in \cite{ramazi2014stability}. 

\begin{proposition}	\label{lem-noEquilibriumIn3DFaces}
	Assume \eqref{snowdriftInequality} holds. 
	The interiors of the faces $\Delta(p^1,p^2,p^4)$, $\Delta(p^1,p^3,p^4)$ and $\Delta(p^2,p^3,p^4)$ do not contain an equilibrium point of the dynamics \eqref{RD}. 
	If $m = 2n+1, n\geq 1$ and $\R = \frac{\T+\S}{2}$, or $m = 2n, n\geq 1$ and $\R = \frac{n\T+(n-1)\S}{2n-1}$, then the interior of the face $\Delta(p^1,p^2,p^3)$ contains the continuum of equilibrium points $\str X^{123}$, and does not contain any other equilibrium.
	For all other values of $m$ and the payoffs, the interior of $\Delta(p^1,p^2,p^3)$ does not contain an equilibrium point.
\end{proposition}

Now we  prove Proposition \ref{lem-equilibria}.

\emph{Proof of Proposition \ref{lem-equilibria}:}
	In view of Proposition \ref{lem-noEquilibriumIn3DFaces}, there is no equilibrium point in the interior of any of $\Delta(p^1,p^2,p^3)$, $\Delta(p^1,p^2,p^4)$, $\Delta(p^1,p^3,p^4)$ and $\Delta(p^2,p^3,p^4)$, except for Cases 3) and 4)  where $\str X^{123}$ appears.
	Hence, all of the rest of the boundary equilibrium points are located on the 6 edges of the simplex. 
	The edges $\Delta(p^1,p^2) = \str X^{12}$ and $\Delta(p^3,p^4) = \str X^{34}$ are always a continuum of equilibrium points. 
	The vertices $p^1,p^2,p^3, p^4$ are also always equilibrium points, but they are included in $\str X^{12}$ and $\str X^{34}$. 
	Hence, the rest of the equilibrium points can be determined by investigating the dynamics in the interior of the remaining four edges, leading to the conclusion (see \cite{hofbauer1998evolutionary}).
\hfill$\blacksquare$

The local stability of the equilibrium points generally depends on the payoffs in $A$, and can be determined based on the convergence results in this section.
However, the following result guarantees the asymptotic stability of $x^{14}$ for all payoffs satisfying \eqref{snowdriftInequality}.
\begin{proposition}		\label{lem-x14x23AreASS}
	Assume \eqref{snowdriftInequality} holds.
	Then $x^{14}$ is asymptotically stable.
\end{proposition}
\begin{IEEEproof}
	The proof follows Proposition \ref{lem-x14x23AreESS} and Lemma \ref{ESSisASS} in Appendix \ref{app:evo}.
\end{IEEEproof}

\subsubsection{Interior equilibrium point}

The dynamics \eqref{RD}, may or may not possess an interior equilibrium depending on the payoff matrix $A$. 
As shown in the following proposition, if the dynamics have an interior equilibrium, it is unique and equal to 
\begin{equation*}
	x^{int} =	
	\begin{bmatrix}
		(a'_{42}-a'_{32})(a'_{13}a'_{24}-a'_{14}a'_{23})	\\
		(a'_{31}-a'_{41})(a'_{13}a'_{24}-a'_{14}a'_{23})	\\
		(a'_{24}-a'_{14})(a'_{31}a'_{42}-a'_{32}a'_{41})	\\
		(a'_{13}-a'_{23})(a'_{31}a'_{42}-a'_{32}a'_{41})
	\end{bmatrix}/r
\end{equation*}
where $a'_{ij}$ are the entries of $A'$ in \eqref{A'}, and 
\begin{align}
	r &= (a'_{13}a'_{24}-a'_{14}a'_{23})(a'_{31}-a'_{41}+a'_{42}-a'_{32}) 	\nonumber\\
		&\ + (a'_{31}a'_{42}-a'_{32}a'_{41})(a'_{13}-a'_{23}+a'_{24}-a'_{14})  > 0	.	\label{r}
\end{align}
The positiveness of $r$ can be derived from \eqref{snowdriftInequality}.
Define the following constants based on the entries $a'_{ij}$ of $A'$:
\begin{gather*}
	b_1 = -\frac{a'_{13}-a'_{23}}{a'_{14}-a'_{24}} = \tfrac{\lceil\frac{m-2}{2}\rceil\S+\lfloor\frac{m}{2}\rfloor\T  - (m-1)\R}{(m-1)(\S-\P)}, \\
	b_2 = -\frac{a'_{42}-a'_{32}}{a'_{41}-a'_{31}} = \tfrac{ \lceil\frac{m-2}{2}\rceil\T+\lfloor\frac{m}{2}\rfloor\S-(m-1)\P}{(m-1)(\T-\R)}. 
\end{gather*}
\begin{proposition}		\label{lem-interior}
	Assume \eqref{snowdriftInequality} holds. 
	It follows that
	\begin{enumerate}
		\item if $\S<\R < \frac{\T+\S}{2}$ or if $m = 2n, n\geq 1$ and $\frac{\T+\S}{2} \leq \R < \frac{n\T+(n-1)\S}{2n-1}$, then the dynamics \eqref{RD} possess exactly one interior equilibrium point $x^{int}$ that is a hyperbolic saddle with two negative eigenvalues; additionally, for all initial conditions on the open line segment
		\begin{equation*}
			\str L^{int} = \left\{ x\in\mathrm{int}(\Delta) \,|\,  x_1 = b_2 x_2, x_4 = b_1 x_3 \right\},
		\end{equation*}		
		the solution trajectory converges to $x^{int}$;
		
		\item otherwise, the dynamics have no interior equilibrium point. 
	\end{enumerate}
\end{proposition}

For the proof, we study the evolution of the ratios $\frac{x_1}{x_2}$ and $\frac{x_4}{x_3}$, which due to the block anti-diagonal structure of the payoff matrix $A'$, are crucial in determining the asymptotic behavior of the replicator dynamics and are explained as follows.

\begin{lemma}		\label{lem-ratios}
	Let $x(0) \in \mathrm{int}(\Delta)$. 
	Then $\frac{d}{dt} \left( \frac{x_1}{x_2}\right)$ is greater than (resp. equal to, resp. less than)  $0$ if and only if 
	$\frac{x_4}{x_3}$ is greater than (resp. equal to, resp. less than) $b_1$.
	Similarly, $\frac{d}{dt} \Big( \frac{x_4}{x_3}\Big)$ is greater than (resp. equal to, resp. less than)  $0$ if and only if 
	$\frac{x_1}{x_2}$ is greater than (resp. equal to, resp. less than) $b_2$.  	
\end{lemma}

\begin{IEEEproof}
	In view of Lemma \ref{faceIsInvariant}, $x(0) \in \mathrm{int}(\Delta)$ implies $x(t)\in\mathrm{int}(\Delta)$ for all $t$.
	Hence, $0<x_{i}(t)<1, i = 1,\ldots,4$, for all $t$.
	So it is possible to define the ratio $\frac{x_i}{x_j}(t), i,j = 1,\ldots,4$ and calculate its time derivative using \cite[Eq. 3.6]{weibull1997} as
	\begin{equation*}		\label{lem-ratio1-1}
		\frac{d}{dt} \Big( \frac{x_i}{x_j}\Big) = [u(p^i,x)-u(p^j,x)] \frac{x_i}{x_j}. 
	\end{equation*}
	Consider the payoff matrix $A'$ and let $i=1, j=2$ and $i=3,j=4$ to obtain the following two equations
	\begin{align}		
		\frac{d}{dt} \Big( \frac{x_1}{x_2}\Big) 
			&= 	[ \underbrace{(a'_{13}-a'_{23})}_{a'_3}x_3 
				+ \underbrace{(a'_{14}-a'_{24})}_{a'_4}x_4 ] 
				\Big(\frac{x_1}{x_2}	\Big)	,		\label{lem-ratio1-2}		\\
		\frac{d}{dt} \Big( \frac{x_4}{x_3}\Big) 
			&= 	[ \underbrace{(a'_{41}-a'_{31})}_{a'_1}x_1
				+\underbrace{(a'_{42}-a'_{32})}_{a'_2}x_2  ] 
				\Big(\frac{x_4}{x_3}\Big). 			\nonumber	
	\end{align}
	In view of Lemma \ref{lem-signStructure}, $a'_1,a'_4>0$. 
	Hence, because of \eqref{lem-ratio1-2}, 
	\begin{equation*}	
		\frac{d}{dt} \Big( \frac{x_1}{x_2}\Big)  > 0
		\Leftrightarrow
		a'_3 x_3 + a'_4 x_4 > 0 
		\xLeftrightarrow[]{a'_4>0}
		\frac{x_4}{x_3} > -\frac{a'_3}{a'_4} = b_1,
	\end{equation*}
	\begin{equation*}	
		\frac{d}{dt} \Big( \frac{x_4}{x_3}\Big)  > 0
		\Leftrightarrow
		a'_1 x_1 + a'_2 x_2 > 0 
		\xLeftrightarrow[]{a'_1>0}
		\frac{x_1}{x_2} > -\frac{a'_2}{a'_1} = b_2.
	\end{equation*}
	This proves the `` greater than" cases. 		
	The `` equal to" and `` less than" cases can be proven similarly.
\end{IEEEproof}

Determining the signs of $b_1$ and $b_2$ will prove useful, and is clarified in the following lemma.
\begin{lemma}	\label{lem-signOfb1b2}
	It holds that $b_2>0$. 
	Moreover, $b_1>0$ (resp. $b_1 = 0$ and $b_1<0$) if and only if  $a'_{13}<a'_{23}$ (resp. $a'_{13}=a'_{23}$ and $a'_{13}>a'_{23}$) where $a'_{ij}$ are the entries of $A'$ in \eqref{A'}.
\end{lemma}
\begin{IEEEproof}
	In view of Lemma \ref{lem-signStructure}, $a'_{32}>a'_{42}$ and $a'_{41}>a'_{31}$. 
	Hence, $b_2>0$ regardless of the payoffs in $A'$.
	Moreover, the inequality $a'_{14}>a'_{24}$ also always holds, which leads to the proof.
\end{IEEEproof}

Now we proceed to the proof of Proposition \ref{lem-interior}.

\emph{Proof of Proposition \ref{lem-interior}: }
	Consider Case 1).
	In view of Lemma \ref{lem-signOfb1b2} and Lemma \ref{lem-signStructure}, $b_1,b_2>0$.
	Then each of the following two sets define a plane in the simplex
	\begin{equation*}		
		\mathcal{P}^{1} = \left\{ x\in\Delta\,|\,	 \frac{x_4}{x_3}  = b_1	\right\},	\quad
		\mathcal{P}^2 = \left\{ x\in\Delta\,|\,	 \frac{x_1}{x_2}  = b_2	\right\}.
	\end{equation*}
	In view of Lemma \ref{lem-ratios}, on each side of the plane $\str P^1$ (resp. $\str P^2$), the quantity $\frac{x_1}{x_2} $ (resp. $\frac{x_4}{x_3} $) either increases or decreases. 
	Hence, if an interior equilibrium point exists, it has to lie on the interior of the intersection of the two planes $\str P^1$ and $\str P^2$, which is the open line segment $\str L^{int}$.
	According to Lemma \ref{lem-ratios}, $\str L^{int}$ is invariant under the replicator dynamics \eqref{RD}.
	The dynamics of $x_2$ on $\str L^{int}$ can be expressed as
	\begin{equation}	\label{lem-interior-1}
		\dot{x}_2 =k  (f x_2 - g) (r x_2 - s) x_2
	\end{equation}
	where
	\begin{align*}
		k &= \frac{1}{(a'_{41}-a'_{31})^2(a'_{13}-a'_{23}+a'_{24}-a'_{14})} > 0,	\\ 
		f &= a'_{32}-a'_{42} + a'_{41}-a'_{31} > 0,		\\
		g &= a'_{41}-a'_{31} > 0,		\
		s = (a'_{13}a'_{24}-a'_{14}a'_{23})(a'_{31}-a'_{41}) > 0,
	\end{align*}
	and $r$ is defined in \eqref{r}.
	The equilibrium points of \eqref{lem-interior-1} are $x^*_2 = 0,\frac{s}{r},\frac{g}{f}$, which are easily proven to be unstable, stable and unstable, respectively.
	On the other hand, $x^*_2 = 0 $ and $x^*_2 = \frac{g}{f}$ correspond to equilibrium points on the boundary of $\Delta$.
	Hence, for any initial condition on $\str L^{int}$, the trajectory $x(t)$ converges to $x^*\in\str L^{int}$ where $x^*_2 = \frac{s}{r}$.
	By using the constraints $\sum_{i=1}^4 x^*_i =1 $ and $x^*\in\str L^{int}$, we get that $x^* = x^{int}$.
	Hence, $x^{int}$ is an interior equilibrium, and for all $x(0)\in\L^{int}$, $x(t)\to x^{int}$. 
	Now the eigenvalues of $x^{int}$ are determined.
	Consider the replicator dynamics \eqref{RD}.
	Replace the vector $x$ by $\hat{x} = \begin{bmatrix} x_1 & x_2 &x_3 &1-x_1-x_2-x_3\end{bmatrix}^\top $, and eliminate the differential equation for $\dot{x}_4$ to get a 3rd order system. 
	Then, the characteristic equation of the corresponding Jacobian matrix about $x^{int}$ is 
	$\lambda^3 + a \lambda^2 + b \lambda + c = 0$ where $a,b,c,\in\mathbb{R}$.
	It can be verified that $c = ab$ and $a > 0 >b,c$.
	Hence, the corresponding eigenvalues of $x^{int}$ are $-a,\pm \sqrt{-b}$, which completes the proof of this case.

	 Now consider Case 2) where $a'_{13} \geq a'_{23}$.
	 Hence, $b_1 \leq 0$ in view of Lemma \ref{lem-signOfb1b2}.
	 Hence, $\str P^1$ does not intersect $\Delta$ implying that the ratio $\frac{x_4}{x_3}$ is always greater than $b_1$.
	 Hence, in view of Lemma \ref{lem-ratios}, $\frac{x_1}{x_2}$ monotonically increases in $\int(\Delta)$. 
	 Hence, there is no interior equilibrium point in this case. 
\hfill$\blacksquare$
\subsection{Trajectories starting on an edge}		\label{subsec-edges}
Due to invariance, the convergence analysis of the dynamics \eqref{RD} on an edge $\Delta(p^k,p^j),k,j\in\{1,2,3,4\}, k\neq j$, can be reduced to the analysis of the following 2-dimensional replicator dynamics
\begin{equation*}
	\dot{x}_i = [(p^i)^\top \hat{A} \hat{x} -\hat{x}^\top \hat{A} \hat{x}]x_i, \qquad i=k,j
\end{equation*}
 where
\begin{equation*}
	\hat{x} = \begin{bmatrix}		x_k	\\ x_j	\end{bmatrix}, \quad 
	\hat{A}_{kj} 
	= \begin{bmatrix} a_{kk}	&a_{kj}	\\	a_{jk}		&a_{kk}	\end{bmatrix}.
 \end{equation*} 
See \cite[Section 3.1.4]{weibull1997},  \cite{hofbauer1998evolutionary} or \cite{ramazi2014stability} for the analysis of these dynamics.

\subsection{Trajectories starting in the interior of a planar face}
We limit this section to the following convergence result that can be proven using the findings in \cite{bomze1, hofbauer1998evolutionary, ramazi2014stability}.
\begin{proposition}	\label{lem-convergenceToPointIn3D}
	If $x(0)$ belongs to one of the faces $\Delta(p^1,p^2,p^3)$, $\Delta(p^1,p^2,p^4)$, $\Delta(p^1,p^3,p^4)$ or $\Delta(p^2,p^3,p^4)$, then $x(t)$ converges to a point in that face as $t\to\infty$.
\end{proposition}

\subsection{Trajectories starting in the interior of the simplex}

\subsubsection{Dynamics in the four sections made by the two ratios} \label{subsec-interiorTrajectories}
When $b_1$ and $b_2$ are positive, the ratios $\frac{x_1}{x_2}$ and $\frac{x_4}{x_3}$ divide the simplex into the four following zones: 
\begin{gather*} 
	\str D^{14} = \left\{ x\in\mathrm{int}(\Delta)\,|\, \frac{x_1}{x_2} > b_2, \ \frac{x_4}{x_3} > b_1		\right\}, \\
	\str D^{23} = \left\{ x\in\mathrm{int}(\Delta)\,|\, \frac{x_4}{x_3}  < b_2 ,\ 	\frac{x_4}{x_3} < b_1 \right\}.	\\
		\str Y^{14} = \left\{ x\in\mathrm{int}(\Delta)\,|\, \frac{x_1}{x_2} > b_2 ,\ 	\frac{x_4}{x_3} < b_1 \right\},	\\
		\str Y^{23} = \left\{ x\in\mathrm{int}(\Delta)\,|\, \frac{x_1}{x_2} < b_2 ,\ 	\frac{x_4}{x_3} > b_1 \right\}.		
\end{gather*}
We investigate interior trajectories of the simplex by studying the dynamics in the above zones and start by $\str D^{14}$ and $\str D^{23}$. 
\begin{lemma}			\label{lem-invariantD}
	$\str D^{14}$ and $\str D^{23}$ are positively invariant under \eqref{RD}.
\end{lemma}
\begin{IEEEproof}	
	First $\str D^{14}$ is shown to be positively invariant.
 	Assume the contrary, i.e., a trajectory $x(t)$ starts from some point in $\str D^{14}$ at $t=t^0$ but does not belong to $\str D^{14}$ at some time $t^*>t^0$. 
	Due to the continuity of the trajectory, there exists some time $t^1 \in (t^0,t^*)$ at which the trajectory intersects the boundary of $\str D^{14}$.
	Hence, at least one of the followings happen
	\begin{equation*}		\label{lem-ratio1-8}
		\frac{x_1}{x_2}(t^1) = b_2, 			\quad		\frac{x_4}{x_3}(t^1) = b_1.
	\end{equation*}
	%Consider the case when only one of the above situations happen.
	Without loss of generality, assume the first case happens.
	Then $\frac{x_1}{x_2}(t^1) < \frac{x_1}{x_2}(t^0) $.
	Hence, $\frac{d}{dt} \Big( \frac{x_1}{x_2}\Big) $ must be negative at some time $t^2\in(t^0,t^1)$. 
	Hence, due to the continuity of the time-derivative of $\frac{x_1}{x_2}$, $\frac{d}{dt} \Big( \frac{x_1}{x_2}\Big) $ is zero at some time $t^3 \in (t^0,t^2)$.
	Hence, in view of Lemma \ref{lem-ratios}, $\frac{x_4}{x_3}(t^3) = b_1$. 	
	This implies that the trajectory has intersected the boundary of $\str D^{14}$ at some time earlier than $t^1$, a contradiction.
	Hence, if a trajectory starts in $\str D^{14}$ at some time $t = t^0$, it remains there afterwards.
	Similarly the positive invariance of $\str D^{23}$ can be shown.
\end{IEEEproof}

\begin{proposition}		\label{lem-convergenceInD}
	Consider a trajectory $x(t)$ of the dynamics \eqref{RD} that passes through $x^0$ at some time $t^0$.
	If $x^0\in \str D^{14}$, then one of the following cases happen
	\begin{equation*}
		\lim_{t\to\infty} x(t) = x^{14} \quad\text{ or } \quad \lim_{t\to\infty} x(t) = x^* \in \str X^{12}\cap \Delta^{NE}.
	\end{equation*}
	If $x^0\in \str D^{23}$, then
	\begin{gather*}
		\lim_{t\to\infty} x(t) = x^* \in (\{x^{23}\}\cup\str X^{12})\cap \Delta^{NE}.
	\end{gather*}	
\end{proposition}
\begin{IEEEproof}
	Consider the case when $x^0\in \str D^{14}$. 
	In view of Lemma \ref{lem-invariantD}, $x(t) \in \str D^{14}$ for all $t\geq t^0$.
	Hence, both inequalities $\frac{x_1}{x_2}(t) > b_2$ and $\frac{x_4}{x_3}(t)> b_1$ hold for all $t\geq t^0$. 				
	Hence, in view of Lemma \ref{lem-ratios}, both ratios $\frac{x_4}{x_3}$ and $\frac{x_1}{x_2}$ monotonically increase with time.
	Hence, each ratio converges to either a constant or $\infty$.
	In case one of the ratios, e.g., $\frac{x_1}{x_2}$, converges to a constant, that constant must be strictly positive.
	This follows from the fact that $\frac{x_1}{x_2}(t^0) > 0$ and that $\frac{x_1}{x_2}$ monotonically increases.
	In general, one of the following cases may occur:

	 1) $\frac{x_1}{x_2}\to \alpha>0$ and $\frac{x_4}{x_3}\to\beta>0$. 
	Thus, $x$ converges to the following line segment
	\begin{equation*}
		\str L^{\alpha,\beta} = \left\{ x\in\Delta \,|\, {x_1} = \alpha {x_2}, {x_4}= \beta {x_3}\right\}.
	\end{equation*}
	In view of Theorem  \ref{lem-wLimitOfLine} in Appendix \ref{app:line}, $x\to\str L^{\alpha,\beta}\cap \Delta^o$.  
	In what follows, it is shown that $\int(\str L^{\alpha,\beta} )\cap\Delta^o=\emptyset$.
	First note that $\alpha>b_2$. 
	This can be proven by contradiction:
	Assume that $\alpha \leq b_2$. 
	Since $x(t^0)\in \str D^{14}$, it holds that $\frac{x_1}{x_2}(t^0)  > b_2$.
	Hence, $\frac{x_1}{x_2}(t^0) > b_2 \geq \alpha$. 
	Then, due to the continuity of the trajectory, there exists some time $t^1>t^0$ such that $\frac{x_1}{x_2}(t^1) = b_2$.
	Hence, $x(t^1)\not\in\str D^{14}$, which contradicts the invariance property of $\str D^{14}$.
	So $\alpha>b_2$.
	%Similarly $\beta>b_1$.		
	Now note that $\int(\L^{\alpha,\beta} ) \subseteq \int(\Delta)$.
	On the other hand, in view of Lemma \ref{lem-interior}, the only interior equilibrium of the system (if there exists any), belongs to the plane $\left\{x\in\Delta\,|\, \frac{x_1}{x_2} = b_2\right\}$. 
	However, as it was discussed above, $\frac{x_1}{x_2}\to\alpha>b_2$.
	Hence, $\int(\Delta)\cap\Delta^o=\emptyset$.
	So $\int(\L^{\alpha,\beta} )\cap\Delta^o=\emptyset$.
	Thus, $x\to\bd(\L^{\alpha,\beta})\cap\Delta^o$.
	The boundary of $\L^{\alpha,\beta}$ consists of the following two points, each of which is an equilibrium:
	\begin{gather*}
		x^{\alpha} = \begin{bmatrix}	\frac{\alpha}{1+\alpha}	&\frac{1}{1+\alpha}	&0	&0	\end{bmatrix}^\top \in\str X^{12}, \\
		 x^{\beta} = \begin{bmatrix} 	0	&0	&\frac{1}{1+\beta}	&\frac{\beta}{1+\beta}	\end{bmatrix}^\top \in\str X^{34}.
	\end{gather*}
	According to Lemma \ref{lem-convergenceToNash} in Appendix \ref{app:Nash}, if $x$ converges to a point, it must belong to $\Delta^{NE}$.
	However, $x^{\beta}\not\in\Delta^{NE}$ in view of Lemma \ref{lem-NashIn3DFaces} in Appendix \ref{app:Nash}.
	Hence, $x\not\to x^{\beta}$ implying that $x\to x^{\alpha}$.
	On the other hand, $x^{\alpha} \in \str X^{12}$ and $x^{\alpha}$ must belong to $\Delta^{NE}$.
	Hence, $x \to x^{*} \in \str X^{12}\cap \Delta^{NE}$.

	2) $\frac{x_1}{x_2}\to \alpha>0$ and $\frac{x_4}{x_3}\to\infty$.
	Hence, $x$ converges to the following line segment
	\begin{equation*}
		\str L^{\alpha,\infty} = \left\{ x\in\Delta \,|\, {x_1}= \alpha{x_2} , {x_3}= 0 \right\}.
	\end{equation*}
	Due to Theorem  \ref{lem-wLimitOfLine}, $x$ converges to an equilibrium or a continuum of equilibria on $\L^{\alpha,\infty}$.  
	On the other hand, $\str L^{\alpha,\infty}$ lies on the face $\Delta(p^1,p^2,p^4)$, and in view of Proposition \ref{lem-noEquilibriumIn3DFaces}, no interior equilibrium exists on this face. 
	Hence, $x$ converges to the intersection of $\str L^{\alpha,\infty} $ with the boundary of $\Delta(p^1,p^2,p^4)$ which is $\{x^{\alpha},p^4\}$.
	However, $p^4\not\in\Delta^{NE}$ and hence $x\not\to p^4$ in view of Lemma \ref{lem-convergenceToNash}.
	Hence, $x\to x^\alpha$.
	So, similar to the previous case, $x \to x^{*} \in \str X^{12}\cap \Delta^{NE}$. 
	
	3) $\frac{x_1}{x_2}\to \infty$ and $\frac{x_4}{x_3}\to\beta>0$.
	Similar to the previous case, it can be shown that $x\to x^{\beta}$ or $x\to p^1$.
	However, neither $x^{\beta}$ nor $p^1$ belongs to $\Delta^{NE}$.
	Hence, this case never happens.
	
	4) $\frac{x_1}{x_2}\to \infty$ and $\frac{x_4}{x_3}\to\infty$.
	Hence, $x$ converges to the following line segment 
	\begin{equation*}
		\str L^{\infty,\infty} = \left\{ x\in\Delta \,|\, {x_2} = 0, {x_3}= 0 \right\} = \Delta(p^1,p^4).
	\end{equation*}
	Due to Theorem  \ref{lem-wLimitOfLine}, $x\to\Delta(p^1,p^4)\cap \Delta^o = \{p^1,x^{14},p^4\}$.
	On the other hand, $p^1,p^4\not\in\Delta^{NE}$.
	Hence, $x\to x^{14}$ in view of Lemma \ref{lem-convergenceToNash}.
	
	Summarizing the above four cases completes the proof for when $x^0\in \str D^{14}$.
	Now let $x^0\in \str D^{23}$. 
	By following the procedure for when $x^0\in \str D^{14}$, it can be shown that both ratios $\frac{x_4}{x_3}$ and $\frac{x_1}{x_2}$ converge either to a positive constant or to $0$.
	In general, one of the following cases may occur:

	 1*) $\frac{x_1}{x_2}\to \alpha>0$ and $\frac{x_4}{x_3}\to\beta>0$. 
	Similar to when $x^0\in \str D^{14}$, this case results in $x \to x^{*} \in \str X^{12}\cap \Delta^{NE}$.

	2*) $\frac{x_1}{x_2}\to \alpha>0$ and $\frac{x_4}{x_3}\to 0$.
	Hence, $x$ converges to the following line segment
	\begin{equation*}
		\str L^{\alpha,0} = \left\{ x\in\Delta \,|\, {x_1} = \alpha{x_2}, {x_4}= 0 \right\}.
	\end{equation*}
	In view of Theorem  \ref{lem-wLimitOfLine}, $x\to\str L^{\alpha,0}\cap \Delta^o$.  
	Clearly $\str L^{\alpha,0} \subseteq \Delta(p^1,p^2,p^3)$.	
	On the other hand, according to Proposition \ref{lem-noEquilibriumIn3DFaces}, $\mathrm{int}(\Delta(p^1,p^2,p^3))\cap\Delta^o$  either is empty or equals to $\str X^{123}$.
	In view of Theorem  \ref{lem-equilibria}, the second case only happens when $m=2n+1,n\geq1$ and $\R=\frac{\T+\S}{2}$, or $m=2n,n\geq1$ and $\R = \frac{n\T+(n-1)\S}{2n-1}$.
	However, for both of these values of $\R$, it can be verified that $b_1<0$.
	Hence, $\str D^{23}=\emptyset$, which contradicts the assumption $x^0\in\str D^{23}$.
	Hence, $\mathrm{int}(\Delta(p^1,p^2,p^3))\cap\Delta^o = \emptyset$.
	So $\int(\L^{\alpha,0}) \cap \Delta^o = \emptyset$ and $x\to\bd(\L^{\alpha,0})$.
	Thus, $x\to\{x^{\alpha},p^3\}$.
	However, $p^3\not\in\Delta^{NE}$ and hence $x\not\to p^3$, in view of Lemma \ref{lem-convergenceToNash}.
	Hence, $x\to x^\alpha$ resulting in $x \to x^{*} \in \str X^{12}\cap \Delta^{NE}$.

	3*) $\frac{x_1}{x_2}\to 0$ and $\frac{x_4}{x_3}\to\beta>0$.
	Hence, $x$ converges to the following line segment
	\begin{equation*}
		\str L^{0,\beta} = \left\{ x\in\Delta \,|\, {x_1}= 0, {x_4}= \beta{x_3} \right\}.
	\end{equation*}
	Similar to the previous case, it can be shown that $x\to \{x^{\beta}, p^2\}$.
	Hence, in view of Lemma \ref{lem-convergenceToNash}, $x\to \{x^{\beta},p^2\} \cap \Delta^{NE}$.
	So $x\to \{p^2\}\cap\Delta^{NE}$ since $x^{\beta}\not\in\Delta^{NE}$.
	On the other hand, $p^2\in\str X^{12}$. 
	Hence, $x\to x^*\in\str X^{12}\cap\Delta^{NE}$.

	4*) $\frac{x_1}{x_2}\to 0$ and $\frac{x_4}{x_3}\to0$.
	Hence, $x$ converges to the following line segment 
	\begin{equation*}
		\str L^{0,0} = \left\{ x\in\Delta \,|\, {x_1} = 0, {x_4}= 0 \right\} = \Delta(p^2,p^3).
	\end{equation*}
	Due to Theorem  \ref{lem-wLimitOfLine}, $x\to\Delta(p^2,p^3)\cap \Delta^o = \{p^2,x^{23},p^3\}$.
	On the other hand, $p^3\not\in\Delta^{NE}$.
	Hence, $x\to \{x^{23},p^2\}\cap \Delta^{NE}$ in view of Lemma \ref{lem-convergenceToNash}.
	Since $p^2\in\str X^{12}$, it can be concluded that $x\to x^*\in (\str X^{12}\cup \{x^{23}\})\cap \Delta^{NE}$.
	
	By summarizing the above cases, the proof for when $x^0\in \str D^{23}$ is complete.
\end{IEEEproof}

\begin{lemma}			\label{lem-convergenceInY}
	Consider a trajectory $x(t)$ of the dynamics \eqref{RD} that passes through $x^0$ at some time $t^0$.
	If $x^0\in\str Y^{14}$, then either $x(t)$ leaves $\str Y^{14}$ after some finite time, or 
	\begin{equation*}
		\lim_{t\to\infty} x(t) = x^{int} \quad\text{or}\quad \lim_{t\to\infty} x(t) = x^* \in \str X^{12}\cap \Delta^{NE}.
	\end{equation*}
	If $x^0\in\str Y^{23}$, then either $x(t)$ leaves $\str Y^{23}$ after some finite time, or 
	\begin{equation*}
		\lim_{t\to\infty} x(t) = x^{int} \ \ \text{or}\ \ \lim_{t\to\infty} x(t) = x^* \in (\str X^{12}\cup\str X^{123})\cap \Delta^{NE}.
	\end{equation*}
\end{lemma}	
\begin{IEEEproof}
	Consider the case when $x^0\in \str Y^{14}$. 
	If $x$ leaves $\str Y^{14}$ after some finite time, the conclusion can be drawn directly.
	So let $\str Y^{14}$ be invariant. 			
	Then the inequalities in the definition of $\str Y^{14}$ hold for all $t\geq t^0$.
	Hence, in view of Lemma \ref{lem-ratios}, $\frac{x_1}{x_2}$ monotonically decreases and hence converges to a constant $\alpha\geq b_2$, and $\frac{x_4}{x_3}$ monotonically increases and hence converges to a constant $\beta\leq b_1$ as $t\to\infty$.
	Hence, $x(t)$ converges to the line segment 
	$
		\L^{\alpha,\beta} = \{ x\in\Delta\,|\, x_1 = \alpha x_2, x_4 = \beta x_3 \}.
	$	
	So based on Theorem  \ref{lem-wLimitOfLine} in Appendix \ref{app:line}, $x(t)$ converges to $\L^{\alpha,\beta}\cap\Delta^o$.
	On the other hand, $\Delta^o$ includes at most one interior equilibrium point $x^{int}$ according to Lemma \ref{lem-wLimitOfLine}.
	Hence, either $x(t)\to x^{int}$ or $x(t)\to\L^{\alpha,\beta}\cap \Delta^{oo}$. 	
	The first case leads to the conclusion directly, so consider the second case.
	First note that $\alpha>0$ since $b_2>0$ in view of Lemmas \ref{lem-signOfb1b2} and \ref{lem-signStructure}.
	Moreover, $\beta>0$ since $\frac{x_4}{x_3}$ monotonically increases from $\frac{x_4}{x_3}(0)>0$ to $\beta$. 	
	Hence, $\alpha, \beta>0$.
	So on the set $\L^{\alpha,\beta}\cap \bd(\Delta)$, either $x_1=x_2 = 0$ or $x_3=x_4=0$ holds.
	Then $\L^{\alpha,\beta}\cap \bd(\Delta) $ equals a point $x^*\in\str X^{12}\cup\str X^{34}$.
	On the other hand, $\Delta^{oo}\subseteq \bd(\Delta)$. 
	Hence, since $\str X^{12}\cup\str X^{34}\subseteq\Delta^{oo}$ it holds that $\L^{\alpha,\beta}\cap \Delta^{oo} = x^*\in\str X^{12}\cup\str X^{34}$.
	Thus, in view of Lemma \ref{lem-convergenceToNash} in Appendix \ref{app:Nash}, $x(t)\to x^*\in(\str X^{12}\cup\str X^{34})\cap\Delta^{NE}$.
	On the other hand, $\str X^{34}\cap\Delta^{NE}=\emptyset$ according to Lemma \ref{lem-NashIn3DFaces} in Appendix \ref{app:Nash}.
	Hence, $x(t)\to x^*\in\str X^{12}\cap\Delta^{NE}$, which completes the proof of this part.

	Now consider the case when $x^0\in \str Y^{23}$ and $\str Y^{23}$ is invariant (otherwise, the result is trivial).
	Hence, in view of Lemma \ref{lem-ratios}, $\frac{x_1}{x_2}$ monotonically increases and hence converges to a constant $\alpha\leq b_2$, and $\frac{x_4}{x_3}$ monotonically decreases and hence converges to a constant $\beta\geq b_1$ as $t\to\infty$.
	So similar to the previous case, either $x(t)\to x^{int}$ or $x(t)\to\L^{\alpha,\beta}\cap \Delta^{oo}$. %xˆ*
	Again the first case leads to the conclusion directly, so consider the second.
	It must be true that $\alpha>0$ since $\frac{x_1}{x_2}$ monotonically increases from $\frac{x_1}{x_2}(0)>0$ to $\alpha$. 	
	If $\beta$ is also positive, then the same as when $x^0\in \str Y^{14}$ takes place, which makes the result trivial.
	So let $\beta =0$.
	Then $\L^{\alpha,\beta}= \{x\in\Delta\,|\, x_1 = \alpha x_2, x_4 = 0 \}$.
	Hence, in view of Theorem  \ref{lem-equilibria}, $\L^{\alpha,\beta}\cap \Delta^{oo} =x^* \in\{ x^{13}\}\cup \str X^{12}\cup \str X^{123}$.
	So in view of Lemma \ref{lem-convergenceToNash} and \ref{lem-singletonEquilibriaNash} in Appendix \ref{app:Nash}, $x(t)\to x^* \in(\str X^{12}\cup \str X^{123})\cap \Delta^{NE}$, which completes the proof.
\end{IEEEproof}

\subsubsection{Global results}
We proceed to the global convergence analysis.
As one would expect, the convergence results depend on the payoffs and to some extent also on $m$.
We provide the results from small to large $\R$ via the following four theorems.
\begin{theorem}	\label{th-smallRConvergence}
	Assume \eqref{snowdriftInequality} holds. 
	Let $x(0)\in\mathrm{int}(\Delta)$.
	Denote the 2-dimensional stable manifold of $x^{int}$ by $W^s(x^{int})$.
	If $\S< \R < \frac{\T+\S}{2} $, then 
	\begin{enumerate}
		\item $x(0)\in W^s(x^{int}) \Rightarrow \displaystyle\lim_{t\to\infty}x(t) = x^{int}$;

		\item $x(0)\not\in W^s(x^{int}) \Rightarrow \displaystyle\lim_{t\to\infty}x(t) = x^{14}$ or $x^{23}$;

		\item $x^{14}$ and $x^{23}$ are asymptotically stable and their basins of attraction are separated by $W^s(x^{int})$;

		\item  $
			x(0) \in \str D^{14} 
			\Rightarrow 
			\displaystyle\lim_{t\to\infty} x(t) = x^{14}
			$ ;
		\item
			$
			 x(0) \in \str D^{23}
			 \Rightarrow
			 \displaystyle\lim_{t\to\infty} x(t) = x^{23}.	
			$
	\end{enumerate}
\end{theorem}		
\begin{IEEEproof}
	Case 1) of the theorem is a direct result of Theorem  \ref{lem-interior}. 
	Now we proceed to Case 2).
	According to Lemma \ref{lem-signOfb1b2} and Lemma \ref{lem-signStructure}, $b_1,b_2>0$.
	Hence, the interior of the simplex can be written as 
	\begin{equation}	\label{th-smallRConvergence-3}
		\mathrm{int}(\Delta) 
		= \str D^{14} \cup \str D^{23} \cup \str Y^{14} \cup \str Y^{23} \cup \str L^{int} \cup \hat{\str P}_{11} \cup \hat{ \str P}_{12} \cup \hat{\str P}_{21} \cup \hat{ \str P}_{22}
	\end{equation}
	where $\L^{int}$ is defined in Theorem  \ref{lem-interior} and
	\begin{align}
		\hat{\str P}_{11} &= \left\{ x \in \int(\Delta) \,|\, \frac{x_4}{x_3} = b_1, \frac{x_1}{x_2} > b_2\right\} ,\label{th-smallRConvergence-4} \\ 
		\hat{\str P}_{12} &= \left\{ x \in \int(\Delta) \,|\, \frac{x_4}{x_3} = b_1, \frac{x_1}{x_2} < b_2\right\}, \nonumber\\
		\hat{\str P}_{21} &= \left\{ x \in \int(\Delta) \,|\, \frac{x_1}{x_2} = b_2,\frac{x_4}{x_3} > b_1\right\}, \nonumber\\
		\hat{\str P}_{22} &= \left\{ x \in \int(\Delta) \,|\, \frac{x_1}{x_2} = b_2,\frac{x_4}{x_3} < b_1\right\}. \nonumber		
	\end{align}
	Hence, $x(0)$ belongs to one of the sets on the right hand side of \eqref{th-smallRConvergence-3}.
	If $x(0)\in \str D^{14}$, then in view of Proposition \ref{lem-convergenceInD}, $x(t)$ converges to either $x^{14}$ or a point in $\str X^{12}\cap \Delta^{NE}$. 
	However, in view of Lemma \ref{lem-NashIn3DFaces} in Appendix \ref{app:Nash},  $\str X^{12}\cap \Delta^{NE}= \emptyset$.
	Hence, $x(t)\to x^{14}$. 
	This proves Case 4).
	Similarly Case 5) can be shown.
	Now consider the case when $x(0)\in \str Y^{14}$. 
	In view of Lemma \ref{lem-convergenceInY}, if $x(t)$ remains in $\str Y^{14}$, it converges to a point in $\str X^{12}\cap \Delta^{NE}$. 
	However, in view of Lemma \ref{lem-NashIn3DFaces},  $\str X^{12}\cap \Delta^{NE}= \emptyset$, which implies $x(t)$ leaves $\str Y^{14}$ after some finite time. 
	Hence, $x(t)$ enters one of the sets $\str L^{int}$, $\hat{\str P}_{11}$, or $\hat{ \str P}_{22}$ at some time $t^1>0$.
	If $x(t^1)\in \str L^{int}$, then $x(t)\to x^{int}$ in view of Proposition \ref{lem-interior}.
	If $x(t^1)\in \hat{\str P}_{11}$, then $x(t)$ enters $\str D^{14}$ after $t=t^1$ since $\frac{x_1}{x_2} > b_2$ in $\hat{\str P}_{11}$ and hence in view of Lemma \ref{lem-ratios}, $\frac{x_4}{x_3}$ increases at $ \hat{\str P}_{11}$.
	So $x(t)\to x^{14}$ in view of Case 4). 
	Similarly, it can be shown that if $x(t^1)\in \hat{\str P}_{22}$, then $x(t)\to x^{23}$.
	Hence, if $x(0)\in \str Y^{14}$, then $x(t)$ converges to one of the points $x^{14},x^{23}$ or $x^{int}$.
	The same can be shown for when $x(0)\in \str Y^{23}$ since $\str X^{123}\not\subseteq\Delta^{oo}$ when $\R<\frac{\T+\S}{2}$. 
	Moreover, the cases when $x(0)$ belongs to one of the sets $\str L^{int} , \hat{\str P}_{11} , \hat{ \str P}_{12} , \hat{\str P}_{21}$ or  $\hat{ \str P}_{22}$ are already included in the arguments for $\str Y^{14}$ and $\str Y^{23}$.
	Hence, $x(t)$ converges to one of $x^{14},x^{23}$ or $x^{int}$.
	On the other hand, only for $x(0)\in W^s(x^{int})$, $x(t)\to x^{int}$. 
	Hence, Case 2) is proven.

	Both $x^{14}$ and $x^{23}$ are asymptotically stable thanks to Proposition \ref{lem-x14x23AreESS} and lemma \ref{ESSisASS}.
	Denote their corresponding basin of attractions by $\str B^{14}$ and $\str B^{23}$.
	Clearly $\str B^{14}$ and $\str B^{23}$ are disjoint.
	Define $\hat{\str B}^{14} := \bd(\str B^{14})\cap\int(\Delta)$ and $\hat{\str B}^{23} := \bd(\str B^{23})\cap\int(\Delta)$.
	Consider a point $x^*\in\hat{\str B}^{14}$.
	The solution $x(t)$ with the initial condition $x^*$, converges to one of $x^{14},x^{23}$ or $x^{int}$ as it was shown above.
	However, $x(t)\not\to x^{14}$ since $x^*\not\in \str B^{14}$.
	Moreover, $x^*\not\in \str B^{23}$ since $x^*\in \bd(\str B^{14})$ and  $\str B^{14}\cap\str B^{23} = \emptyset$ and $\str B^{23}$ is open.
	Hence, $x(t)\to x^{int}$.
	So $x^{int}$ lies on $\str B^{14}$.
	The same can be shown for $\str B^{23}$.
	Now both $\hat{\str B}^{14}$ and $\hat{\str B}^{23}$ are 2-dimensional invariant manifolds, and for any initial condition located on them, $x(t)\to x^{int}$.
	On the other hand, $x^{int}$ is hyperbolic in view of Theorem  \ref{lem-interior}, and hence $W^s(x^{int})$ is the unique 2-dimensional invariant manifold passing through $x^{int}$.
	Hence, $\hat{\str B}^{14}$ and $\hat{\str B}^{23}$ coincide and are equivalent to $W^s(x^{int})$. 
	This proves Case 3) and hence the whole.	
\end{IEEEproof}

\begin{figure}[h]
	\centering
	\includegraphics[width = .5\textwidth]{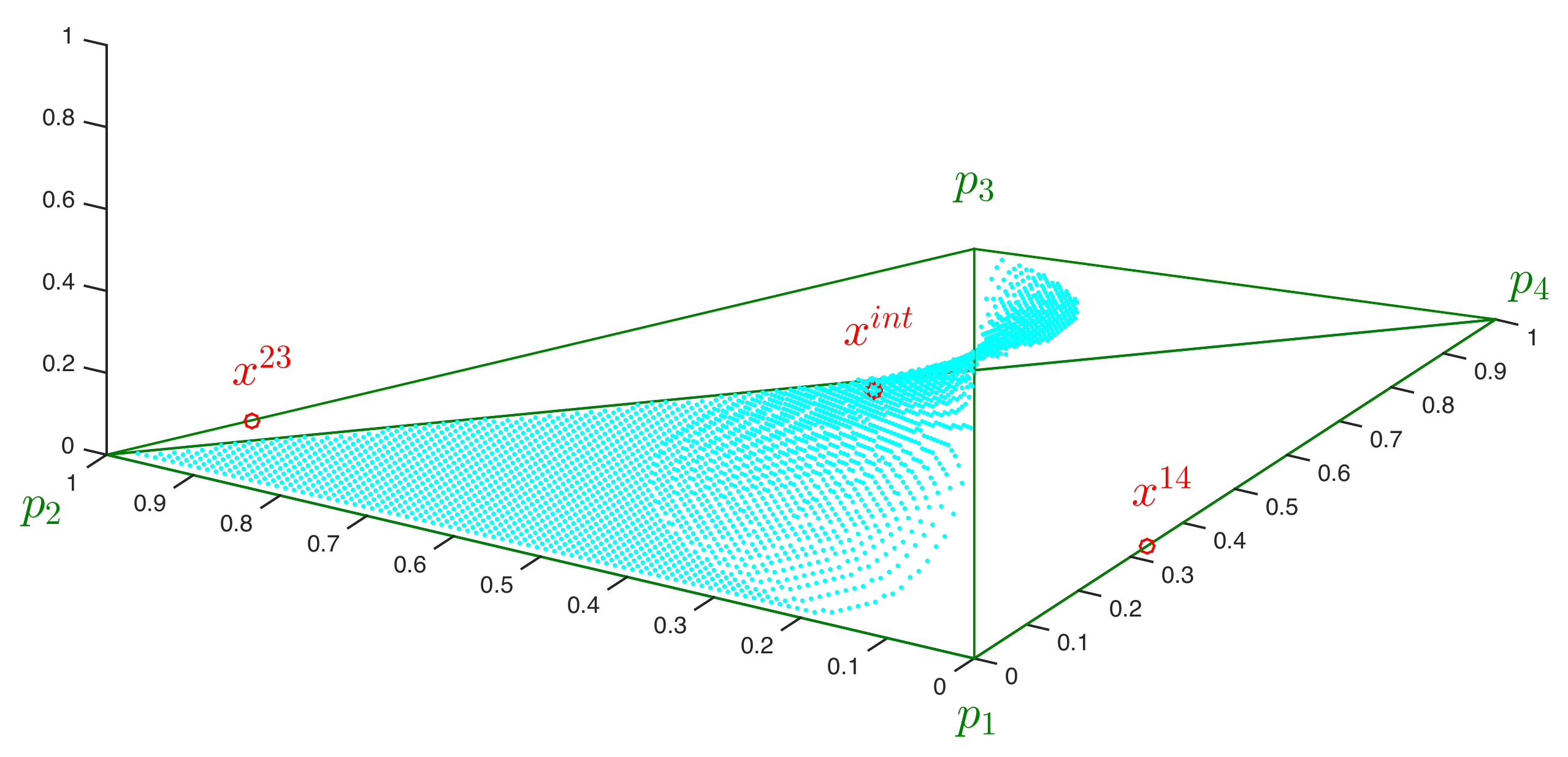}
	\caption{An example of the two-dimensional stable manifold mentioned in Theorem \ref{th-smallRConvergence} for payoff values $\T = 6$, $\R = 4$, $\S = 3$, $\P = 2$ and the number of repetitions $m = 8$. The \textcolor{cyan}{cyan} points are samples of the stable manifold $W^s(x^{int})$.}
	\label{Fig1}
\end{figure}

An example of the two-dimensional stable manifold mentioned in Theorem \ref{th-smallRConvergence} is shown in Figure \ref{Fig1}.
For intermediate values of $\R$, the convergence results depend on whether $m$ is odd or even. 
Therefore, two separate theorems are dedicated to these values.		
\begin{theorem}	\label{th-evenMConvergence}	
	Assume \eqref{snowdriftInequality} holds. 
	Let $x(0)\in\mathrm{int}(\Delta)$.
	Assume $m=2n,n\geq1$.
	It follows that 
	\begin{enumerate}		
		\item if $\R = \frac{\T+\S}{2}$, then 
		\begin{equation*}
			\lim_{t\to\infty} x(t) = x^* \in \{x^{14},x^{int},p^2\};
		\end{equation*}
		
		\item if  $\frac{\T+\S}{2} < \R <  \frac{n\T+(n-1)\S}{2n-1}$, then 	
		\begin{equation*}
			\lim_{t\to\infty} x(t) = x^* \in \{x^{14},x^{int}\} \cup(\str X^{12}\cap\Delta^{NE});
		\end{equation*}
	
		\item if $\R = \frac{n\T+(n-1)\S}{2n-1}$, then 
		\begin{equation*}
			\lim_{t\to\infty} x(t) = x^* \in \{x^{14}\}\cup\str X^{123}\cup(\str X^{12}\cap \Delta^{NE}).
		\end{equation*}
	\end{enumerate}	
\end{theorem}		
\begin{IEEEproof}
	In view of Lemma \ref{lem-signOfb1b2}, $b_1,b_2>0$ in Cases 1) and 2).
	Hence, by following the same steps as in the proof of Theorem  \ref{th-smallRConvergence}, it can be shown that
	$x(t)\to x^*\in\{x^{14},x^{23},x^{int}\}\cup(\str X^{12}\cap \Delta^{NE})$.
	However, $x^{23}\not\in\Delta^{o}$ in view of Theorem  \ref{lem-equilibria} and hence $x(t)\not\to x^{23}$.
	Then in view of Lemma \ref{lem-convergenceToNash}, Case 2) is proven.
	Moreover, the fact that $\str X^{12} \cap \Delta^{NE} = \{p^2\}$ for $\R = \frac{\T+\S}{2}$, proves Case 1). 
	For Case 3), $b_2>0$, but $b_1=0$ in view of Lemma \ref{lem-signOfb1b2}.
	Hence, $\int(\Delta)$ can be written as follows
	\begin{equation*}	\label{th-evenMConvergence-1}
		\mathrm{int}(\Delta) 
		= \str D^{14} \cup \str Y^{14} \cup \hat{\str P}_{11}
	\end{equation*}
	where $\hat{\str P}_{11}$ is defined in \eqref{th-smallRConvergence-4}.
	Then similar to the proof of Theorem  \ref{th-smallRConvergence}, we arrive at the conclusion.
\end{IEEEproof}

\begin{theorem}		\label{th-oddMConvergence}
	Assume \eqref{snowdriftInequality} holds. 
	Let $x(0)\in\mathrm{int}(\Delta)$.
	Assume $m=2n+1,n\geq1$.
	It follows that  
	\begin{enumerate}		
		\item if $\R = \frac{\T+\S}{2}$, then 
		\begin{equation*}
			\lim_{t\to\infty} x(t) = x^* \in \{x^{14},x^{23}\}\cup\str X^{123};
		\end{equation*}
		
		\item if  $\frac{\T+\S}{2} < \R \leq \frac{(n+1)\T+n\S}{2n+1}$, then 	
		\begin{equation*}
			\lim_{t\to\infty} x(t) = x^{14};
		\end{equation*}
	\end{enumerate}	
\end{theorem}	
\begin{IEEEproof}
	In view of Lemma \ref{lem-signOfb1b2}, $b_2>0\geq b_1$ in all cases.
	Hence, by following the same steps as in the proof of Case 3) of Theorem  \ref{th-evenMConvergence}, it can be shown that
	$x(t)\to x^*\in\{x^{14},x^{23}\}\cup(\str X^{12}\cap \Delta^{NE})\cup \str X^{123}$ where $\str X^{123}$ shows up only in Case 1) according to Proposition  \ref{lem-equilibria}.
	Then according to Lemma \ref{lem-NashIn3DFaces} in Appendix \ref{app:Nash}, $\str X^{12}\cap\Delta^{NE} = \emptyset $, which proves Case 1) and Case 2) except for when $\R$ equals $ \frac{(n+1)\T+n\S}{2n+1}$.
	When the equality happens, $\str X^{12}\cap\Delta^{NE} = \{p^2\} $ in view of Lemma \ref{lem-singletonEquilibriaNash} in Appendix \ref{app:Nash}.
	However, in view of Lemma \ref{lem-ratios}, $b_1 \leq 0$ implies that $\frac{x_1}{x_2}$ monotonically increases.
	Hence, $\frac{x_1}{x_2}(t)>\frac{x_1}{x_2}(0)$ for all $t>0$.
	On the other hand, $\frac{x_1}{x_2}(0)>p^2$ since $x(0)\in\int(\Delta)$.
	Hence, $x(t)\not\to p^2$, which completes the proof.
\end{IEEEproof}	

\begin{theorem}	\label{th-largeRConvergence}
	Assume \eqref{snowdriftInequality} holds. 
	Let $x(0)\in\mathrm{int}(\Delta)$.
	If $\max\left\{ \frac{\lceil\frac{m-2}{2}\rceil\S+\lfloor\frac{m}{2}\rfloor\T}{m-1} ,\frac{\lceil\frac{m}{2}\rceil\T+\lfloor\frac{m}{2}\rfloor\S}{m}    \right\}<\R<\T$, then 	
	\begin{equation*}
		\lim_{t\to\infty} x(t) = x^* \in \{x^{14}\} \cup(\str X^{12}\cap\Delta^{NE}).
	\end{equation*}
\end{theorem}			
\begin{IEEEproof}
	The proof is similar to that of Case 2) in Theorem  \ref{th-oddMConvergence}.
\end{IEEEproof}
	
Note that each equilibrium on $\str X^{34}$ performs as an $\alpha$-limit point in the case of Theorem \ref{th-largeRConvergence}.
The integration of the convergence results when the initial condition is in the interior of the simplex and when it is on the boundary of the simplex, yields the following corollary.
\begin{corollary}	\label{th-finalTh}
	Assume \eqref{snowdriftInequality} holds. 
	For any initial condition $x(0)\in\Delta$, the solution $x(t)$ of the replicator dynamics \eqref{RD}, converges to a point in $\Delta$ as time goes to infinity.
\end{corollary}

Therefore, no limit cycle or strange attractor can take place in the dynamics, and we always have convergence to a point. 

\subsection{Discussion}
%%
%\subsection{Evolution of the population ratios of the four types of players}
Now that we know the asymptotic behavior of the replicator dynamics \eqref{RD} for all range of payoffs, we can proceed to the interpretation of the results in terms of the individuals playing the four types of strategies. 
We use two performance measures to compare the population at different states $x$.
The first one is \emph{average population payoff}  $x^\top A x$.
The second is average number of times cooperation is played in the population, which we call the \emph{average cooperation level} and denote by
\begin{equation*}
	x_C := \sum_{i,j\in\{1,\ldots,4\}} x_i x_j \frac{C_{ij} }{2m}
\end{equation*}
where $C_{ij}$ is the number of times cooperation is played in the $m$ rounds when two individuals playing the strategies corresponding to indices $i$ and $j$ are matched to play the repeated game $G^m$.
As an illustration, $C_{11} = 2m$, as both matched $ALLC$ players cooperate in every $m$ rounds, and $C_{14} = C_{41} = m$ as only the $ALLC$ player cooperates when matched with an $ALLD$ player.
Moreover, the average cooperation level at $x^{14}$ is $\frac{\S-\P}{\S-\P+\T-\R}$ since only $ALLC$ players cooperate, and reaches 1 at any state in $\str X^{12}$ since both $ALLC$ and $TFT$ players cooperate.

\begin{figure}[h]
	\centering
	\includegraphics[width = .45\textwidth]{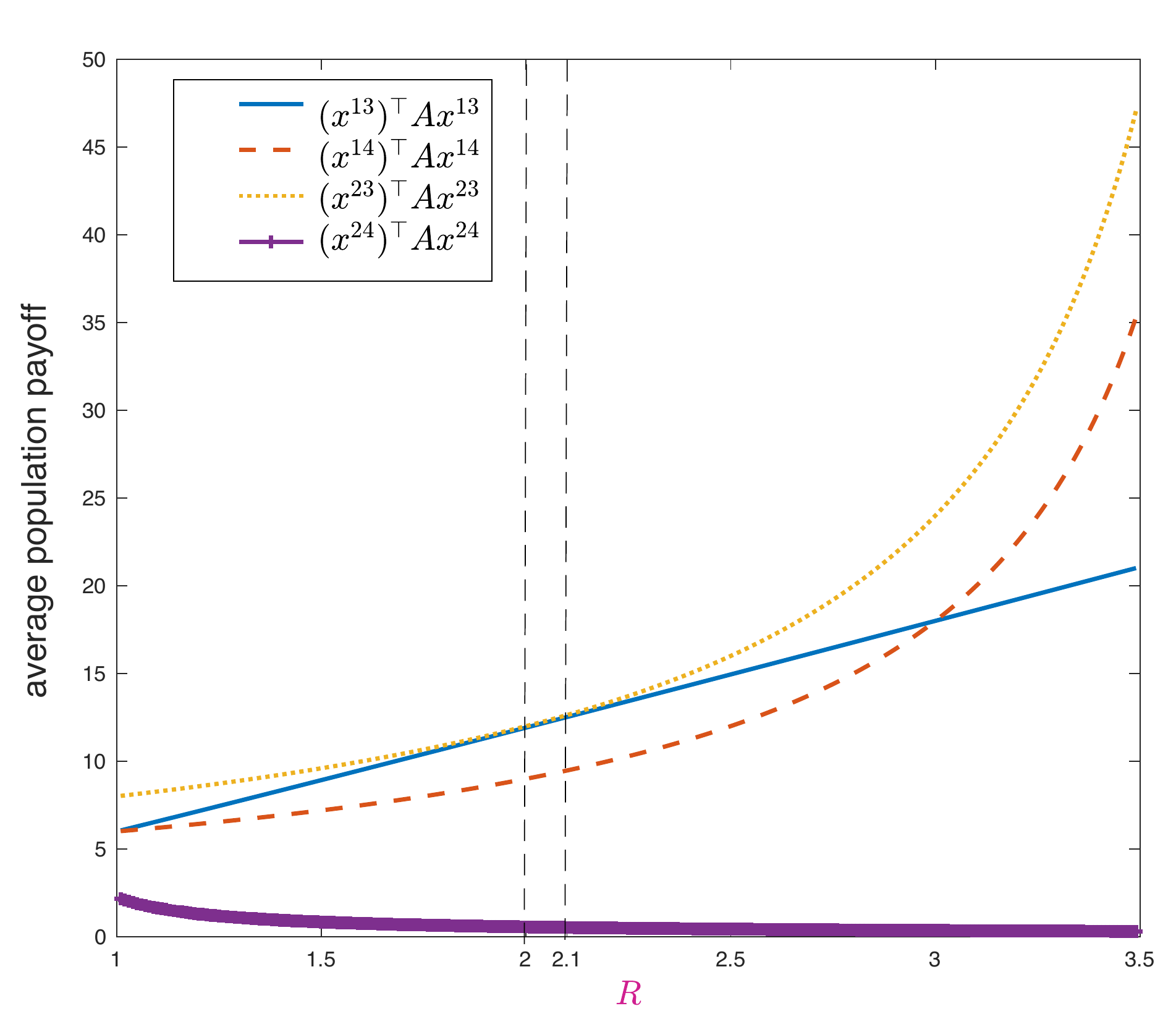}
	\caption{\textbf{Average population payoff at the equilibria as a function of $\R$.}
	Other parameters are set to be $\T = 3, \S = 1, \P = 0$ and $m = 6$.}	\label{Fig2}
\end{figure}	
%
%Now assume that the payoffs satisfy the snowdrift constraint \eqref{snowdriftInequality}.
Now consider a population where the portions of individuals playing $ALLC$, $TFT$, $STFT$ and $ALLD$ are all nonzero.
For small values of $\R$, i.e., less than the average of $\T$ and $\S$, almost always the population converges to one of the following states: \emph{i)} $x^{14}$ that is a mixed population of $ALLC$ and $ALLD$ players or 
\emph{ii)} $x^{23}$ that is a mixed population of $TFT$ and $STFT$ players.
Both states are evolutionary (and hence asymptotically) stable (see Appendix-B). 
Therefore, evolutionary forces select against any mutant population at these two states.
Moreover, for a zero-measure set of initial states, the population converges to $x^{int}$ where all four types of players are present. 
However, $x^{int}$ is unstable and small perturbations can lead the population to one of $x^{14}$ and $x^{23}$. 
Between the two, $x^{23}$ has a higher average population payoff since 
\begin{align*}
	(x^{23})^\top Ax^{23} &- (x^{14})^\top Ax^{14}  \\
	&= \begin{cases}
		\frac{m(\S-\T)^2}{4(\T-\R + \S-\P)}>0 & m =2n, n\geq 1 \\
		\frac{(m^2-1)(\S-\T)^2}{4m(\T-\R + \S-\P)}>0 & m = 2n+1, n\geq 1
	\end{cases}
\end{align*}
(see Figure \ref{Fig2}) as well as a higher cooperation level since 
\[
	x^{23}_C - x^{14}_C
	= \begin{cases}
		\frac{\T-\S}{2(\T-\R + \S-\P)}>0 & m =2n, n\geq 1 \\
		\frac{(m-1)(\T-\S)}{2m(\T-\R + \S-\P)}>0 & m = 2n+1, n\geq 1
	\end{cases}.
\]

Now if the base game is repeated for even number of times, as $\R$ increases, the state $x^{23}$ moves towards $p^2$ where only $TFT$ players are present. 
When $\R$ equals the average of $\T$ and $\S$, $x^{23}$ coincides with $p^2$, and hence $STFT$ players stand out (except for those zero-measure initial conditions that lead to $x^{int}$).  
As $\R$ further increases, the single equilibrium state $p^2$ is expanded to the set of a continuum of equilibria $\str X^{12}\cap\Delta^{NE}$.
Therefore, the population either converges to $x^{14}$ where $ALLC$ and $ALLD$ players coexist or to a state where $ALLC$ and $TFT$ players coexist. 
As one would expect, any equilibrium $x^{\alpha}\in\str X^{12}\cap\Delta^{NE}$ outperforms $x^{14}$ in terms of both average population payoff and average cooperation level as 
\[
	(x^{\alpha})^\top Ax^{\alpha} - (x^{14})^\top Ax^{14}
	= \frac{m(\T-\R)(\R-\S)}{(\T-\R + \S-\P)}>0,
\]
and that $x^{\alpha}$ has the highest possible average cooperation level $x^{\alpha}_C = 1$. 

At the same time, $x^{int}$ is moving towards the face $\Delta(p^1,p^2,p^3)$, and when $\R$ equals $\frac{n\T+(n-1)\S}{2n-1}$, $x^{int}$ lies on $\str X^{123}$ where $ALLC$, $TFT$ and $STFT$ players coexist.

If the base game is repeated for odd number of times, $STFT$ players survive for a greater range of $\R$. 
This time for $\R$ being equal to $\frac{\T+\S}{2}$, $x^{int}$ lies on $\str X^{123}$.
Then suddenly, by a small increment in $\R$, the  set $\str X^{123}$ disappears, and no population converges to $x^{23}$. 
Therefore, starting from any initial condition, the population converges to the polymorphic population of $ALLC$ and $ALLD$ players $x^{14}$. 
This is the only situation where both conditional strategies $TFT$ and $STFT$ are wiped out of the population, and is continued up to when $\R$ equals $\frac{(n+1)\T + n\S}{2n+1}$. 
It can be verified that both the average population payoff and cooperation level at $x^{14}$ monotonically increase in $\R$. 
Therefore, as expected, increasing $\R$ results in a more profitable and cooperative long-term population.

When $\R$ further increases, the behavior of the system is almost the same for both odd and even $m$.
The population either converges to $x^{14}$ where  $ALLC$ and $ALLD$ players coexist or to an equilibrium on $\str X^{12}\cap\Delta^{NE}$ where $ALLC$ and $TFT$ players coexist.
So the suspiciousness of $STFT$ players wipes them out from the population.
Moreover, as $\R$ increases, $x^{14}$ gets closer to $p^1$ where all individuals are $ALLC$ players.

In general, perhaps $STFT$ can be considered as the worst strategy in terms of survival especially for $\R>\frac{\T+\S}{2}$.
Conversely, regardless of the payoffs, there always exists a set of initial conditions for which $ALLC$ players show up in the long run.
Moreover, in addition to $x^{14}$, all the limit states in $\str X^{12}\cap\Delta^{NE}$ (except for $p^2$) have a nonzero portion of $ALLC$ players.
This, surprisingly, makes the simple unconditional $ALLC$ strategy perhaps the most robust in terms of survival and appearance in the long run. 
This may explain the existence of individuals who unconditionally cooperate in real-life scenarios that can be captured by repeated snowdrift games.

Interestingly, $x^{14}$ is always an evolutionary (and asymptotically) stable state of the system, regardless of the payoffs. 
This state consists of $\frac{\S-\P}{\S-\P+\T-\R}$ $ALLC$ players that can be considered as cooperators and $\frac{\T-\R}{\S-\P+\T-\R}$ $ALLD$ players that can be considered as defectors.
On the other hand, the unique evolutionarily stable state of the base game consists of $\frac{\S-\P}{\S-\P+\T-\R}$ $C$ players, i.e., cooperators, and $\frac{\T-\R}{\S-\P+\T-\R}$ $D$ players, i.e., defectors.
Thus, the repetition of the base game and the introduction of the two conditional strategies $TFT$ and $STFT$, does not eliminate or even change this evolutionarily stable mixture of cooperators and defectors, but adds some new more-cooperative final states such as those on $\str X^{12}$. 

Moreover, since $x^{14}_C = \frac{\S-\P}{\S-\P+\T-\R}$ and $x^{\alpha}_C = 1$ for any $x^{\alpha}\in\str X^{12}$,
adding enough $TFT$ players to a population of $ALLC$ and $ALLD$ players can dramatically increase the average level of cooperation, if $\R$ is large enough.
The claim does not change when $STFT$ players are also present in the population. 
More specifically, if $\R$ is greater than the lower bound provided in Theorem \ref{th-largeRConvergence}, and $x_2(0)$, the initial portion of $TFT$ players is large enough so that $x(0)$ belongs to the basin of attraction of $\str X^{12}$, then the population state converges to a point on $\str X^{12}$ that has a higher average population payoff and cooperation level.

The convergence analysis also reveals how the average cooperation level changes as $\R$ increases.
Particularly, in the presence of the four types of players, increments in $\R$ make the final population more probable to become completely cooperative.

%%%
\section{Concluding remarks}		\label{secConclusion}
	
	Our analysis highlights repetition as a mechanism that promotes cooperation among selfish individuals in snowdrift social dilemmas. 
	Unlike the trend of research on repeated games that allows for a wide range of complicated and uncommon reactive strategies, we have limited them to four typical ones.
	This provides a more realistic setup for human societies \cite{samuelson2003evolutionary}. 
	On the other hand, we have modelled the evolutions of the players' population portions by replicator dynamics which well approximate the behavior of well-mixed large populations governed by the proportional imitation update rule \cite{sandholm2010population}.  
	%Moreover, the fact that the base game is a two strategy snowdrift, implies that the population consists of binary decision-making individuals who do the opposite of what the average population does.
	Given all this, we show that for large well-mixed populations of imitative individuals who play snowdrift games, repeating the game and the introduction of the conditional strategy $TFT$ promotes cooperation.  
	However, this is not because $TFT$ players are long-term dominants as often reported in repeated prisoner's dilemma, but because they lead to more-cooperative final population states which are also more profitable. 
	This promotion of cooperation is preserved even if some of the $TFT$ players start their interactions suspiciously and defect initially; that is, if there are also some $STFT$ players in the population. 
	Indeed, for low values of reward $\R$, such players survive, yet for high rewards they become extinct. 
	Finally, those who always cooperate regardless of their opponents' moves have a high chance of survival, which may explain the observation of such behaviors in real life.

	Our main technical contribution is the study of the ratios of the state variables of the replicator dynamics and to show their monotonicity over time. 
	Such a technique has the potential to be applied to other replicator dynamics whose payoff matrix has repeated entries.   
	We, finally, emphasize that all payoffs and even the number of repetitions are parameters, which not only gives rise to the above outcomes, but more importantly provides a parametric framework to control the final average population payoff and cooperation level.
	This simply can be done by tuning the parameters according to the convergence results in Theorems \ref{th-evenMConvergence} to \ref{th-largeRConvergence}.

%%%
\appendix
\subsection{Lemma \ref{lem-RDInvariantUnderAddition}}	\label{app:A}
\begin{lemma}		\label{lem-RDInvariantUnderAddition}
	The replicator dynamics \eqref{RD} are invariant under the addition of a constant to all of the entries of a column of the payoff matrix $A$.
\end{lemma}
\begin{IEEEproof}
    	See \cite[Section 3.1.2]{weibull1997}.
%	Let $A'$ denote the resulted matrix after the addition of a constant $b\in\mathbb{R}$ to the $i$th column of $A$, i.e.,
%	\[
%		A' =  A + \mathbf{1}p^i b
%	\]
%	where $\mathbf{1}$ is the all-one column vector. Let $v$ be the utility function under $A'$. Then
%	\begin{align*}
%		v(p^i,x) - v(x,x) 			
%		&= p^i^\top  A' x - x^\top  A' x 		\\
%		&= p^i^\top  Ax + p^i^\top  (\mathbf{1}p^i b)x - x^\top  (\mathbf{1}p^i b) x - x^\top  A x			\\
%		&= p^i^\top  Ax + p^i bx - p^i b x - x^\top  A x			\\
%		&= p^i^\top  Ax  - x^\top  A x			\\
%		&= u(p^i,x) - u(x,x).
%	\end{align*}
%	Hence, according to \eqref{RD}, the dynamics remain unchanged with $A'$ in place of $A$. 
\end{IEEEproof}

\subsection{Evolutionary Stability: Proof of Proposition \ref{lem-x14x23AreASS}} \label{app:evo}
A state $x\in\Delta$ is said to be an \emph{evolutionarily stable state (strategy) (ESS) of} $A$ if it satisfies the following two conditions \cite[pp. 81]{sandholm2010population}:
\begin{align}
	&x^\top Ax \geq y^\top Ax		\quad \forall y\in\Delta,							\label{ESSCondition1}	\\
	&[x^\top Ax = y^\top Ax 		\quad\text{and}\quad y\neq x ]\Rightarrow x^\top Ay > y^\top Ay 	\label{ESSCondition2}.
\end{align}
The set of all evolutionarily stable states is denoted by $\Delta^{ESS}$.
\begin{lemma}[Proposition 3.10 in \cite{weibull1997}]		\label{ESSisASS}
	Every $x\in\Delta^{ESS}$ is asymptotically stable under the replicator dynamics \eqref{RD}.
\end{lemma}
\begin{proposition}		\label{lem-x14x23AreESS}
	$x^{14}\in\Delta^{ESS}$.
	Moreover, $x^{23}\in\Delta^{ESS}$ if $\R<\frac{\T+\S}{2}$.
\end{proposition}
\begin{IEEEproof}
	The result for $x^{14}$ is proven in the following, and that for $x^{23}$ can be done similarly. 
	Consider
	\begin{equation*}
		Ax^{14} 
		= \frac{1}{{a'_{14}}+{a'_{41}}}
		\begin{bmatrix}
			{a'_{14}}{a'_{41}}	&	{a'_{24}}{a'_{41}}	&	{a'_{14}}{a'_{31}}	&	{a'_{14}}{a'_{41}}
		\end{bmatrix}^\top .
	\end{equation*}
	In view of Lemma \ref{lem-signStructure}, $	{a'_{41}}>{a'_{31}}\geq 0$ and ${a'_{14}}>{a'_{24}}\geq 0$.
	Hence, ${a'_{14}}{a'_{41}} > {a'_{24}}{a'_{41}}, {a'_{14}}{a'_{41}}$, implying that the maximum element of $Ax^{14}$ is ${a'_{14}}{a'_{41}}$.
	Hence, any $y\in\Delta$ satisfying $y_2,y_3 = 0$ maximizes $y^\top  A x^{14}$.
	So ${x^{14}}^\top  A x^{14}$ is the maximum of $y^\top  A x^{14}$, which implies that \eqref{ESSCondition1} is in force. 
	On the other hand, if for some $y\in\Delta$, ${x^{14}}^\top Ax^{14} = y^\top  A x^{14}$, then $y$ maximizes $y^\top  A x^{14}$.
	Such a $y$ satisfies $y_2,y_3 = 0$, which results in
	\begin{equation}	\label{lem-y1y2AreESS-1}
		y^\top Ax^{14}  = \frac{{a'_{14}}^2y_4 + {a'_{41}}^2 y_1}{{a'_{14}}+{a'_{41}}} ,	\quad
		y^\top Ax^{14}  = ({a'_{14}}+{a'_{41}}) y_1y_4.
	\end{equation}	
	On the other hand,
	\begin{align*}
		&[y_4(a'_{14}+{a'_{41}}) -{a'_{41}}]^2 \geq 0	\\	
		\iff
		&{a'_{14}}^2 y_4 + {a'_{41}}^2 (1-y_4) \geq ({a'_{14}}+{a'_{41}})^2 y_4(1-y_4) 		\\
%		\iff
%		&{a'_{14}}^2 y_4 + {a'_{41}}^2 y_1 \geq ({a'_{14}}+{a'_{41}})^2 y_1y_4		\\
		\iff
		&\frac{{a'_{14}}^2y_4 + {a'_{41}}^2 y_1}{{a'_{14}}+{a'_{41}}} \geq ({a'_{14}}+{a'_{41}}) y_1y_4.
	\end{align*}
	Hence, in view of \eqref{lem-y1y2AreESS-1}, ${x^{14}}^\top Ay \geq y^\top Ay$.
	However, the equality holds only when
	\begin{gather*}
		[y_4({a'_{14}}+{a'_{41}}) -{a'_{41}}]^2 = 0	\\
		\Rightarrow
		y_4 = \frac{{a'_{41}}}{{a'_{14}}+{a'_{41}}}, y_1 = \frac{{a'_{14}}}{{a'_{14}}+{a'_{41}}} 
		\Rightarrow
		y = x^{14}.
	\end{gather*}
	Hence, ${x^{14}}^\top Ax^{14} > y^\top Ax^{14}$ for all $y\neq x^{14}$.
	%Hence, it is shown that when ${x^{14}}^\top A{x^{14}} = y^\top A{x^{14}}$ and $x\neq x^{14}$, it holds that ${x^{14}}^\top Ax > y^\top Ax$.
	So \eqref{ESSCondition2} is true, implying $x^{14}\in\Delta^{ESS}$. 
\end{IEEEproof}

\subsection{Nash equilibria and their relation to convergence points}	\label{app:Nash}

Call a trajectory $x(t)$ an \emph{interior trajectory}, if $x(0)\in \mathrm{int}(\Delta)$.
When investigating the final state (convergence point) of an interior trajectory, several equilibrium points often show up as possible candidates. 
In what follows, a known game theoretical result is reviewed to confine the possible candidates.
Define $\Delta^{NE}$, the subset of strategies (states) that are in Nash equilibrium with themselves \cite[Section 1.5.2]{weibull1997}, by
\begin{equation*}
	\Delta^{NE} = \left\{ x\in\Delta\,|\, x^\top  A x \geq y^\top  A x \quad \forall y\in \Delta \right\} .
\end{equation*}
\begin{lemma}	\label{lem-convergenceToNash} \emph{(\cite[Proposition 3.5]{weibull1997})}
	If an interior trajectory $x(t)$ converges to a point $x^*$, then $x^*\in\Delta^{NE}$. 
\end{lemma}

Similar to Lemma \ref{lem-RDInvariantUnderAddition}, it can be easily verified that $\Delta^{NE}$ is invariant under the addition of a constant to all of the entries of a column of the payoff matrix $A$.
Hence, we change $A$ in the definition of $\Delta^{NE}$ with the more simple-structure payoff matrix $A'$ in future derivations.
The following lemma reveals those points of $\str X^{12}$ and $\str X^{34}$ that belong to $\Delta^{NE}$.
\begin{lemma}	\label{lem-NashIn3DFaces}
	Assume \eqref{snowdriftInequality} holds. 
	Then $\str X^{34}\cap\Delta^{NE} = \emptyset$. 
	Moreover,
	\begin{itemize}
		\item if $\S< \R < \frac{\T+\S}{2} $ or $m=2n+1,n\geq 1$ and $\frac{\T+\S}{2} \leq \R < \frac{(n+1)\T+n\S}{2n+1}$, then $$ \str X^{12} \cap \Delta^{NE} = \emptyset ;$$

		\item if $m=2n+1,n\geq 1$ and $\R = \frac{(n+1)\T+n\S}{2n+1}$, or $m=2n,n\geq 1$ and $\R = \frac{\T+\S}{2}$, then 
			$$ \str X^{12} \cap \Delta^{NE} = \{ p^2\};$$
						
		\item if $m=2n+1,n\geq 1$ and $ \frac{(n+1)\T+n\S}{2n+1}< \R<\T$, or $m=2n,n\geq 1$ and $\frac{\T+\S}{2} < \R < \T$, then 	
			\begin{align*}
				 &\str X^{12} \cap \Delta^{NE} = \\
				 &\scalebox{0.75}{$\left\{ \alpha p^1 + (1-\alpha)p^2\,\Big |\, \alpha \in \left [0, 
			\min\left\{ \tfrac{m\R- \lceil\frac{m}{2}\rceil\T-\lfloor\frac{m}{2}\rfloor\S}{\T-\R}, \tfrac{m\R-\T-(m-1)\P}{m(\T-\R)},1 \right \}\right]\right\}$}.
			\end{align*}
	\end{itemize}
\end{lemma}
\begin{IEEEproof}	
	Let $x\in \str X^{34}$.
	Then 
	$
		A'x = 
		\begin{bmatrix}
			 a'_{13}x_3 + a'_{14}x_4 	
			&a'_{23}x_3 + a'_{24} x_4	
			&0
			&0
		\end{bmatrix}^\top 
	$.
	So based on the definition of $\Delta^{NE}$,
	\begin{equation*}
		x \in \Delta^{NE}
		\iff
		a'_{13}x_3 + a'_{14}x_4 \leq 0 \text{ and } a'_{23}x_3 + a'_{24} x_4\leq 0.
	\end{equation*}
	However, in view of Lemma \ref{lem-signStructure}, $a'_{13},a'_{14},a'_{23},a'_{24}>0$. 	
	Hence, because of $x_3+x_4 =1$ and $x_3,x_4\geq 0$, it can be concluded that $x\not\in \Delta^{NE}$.
	Now let $x\in \str X^{12}$.
	Then 
	\begin{equation*}
		A'x = 
		\begin{bmatrix}
			0 	
			&0	
			&a'_{31}x_1 + a'_{32}x_2 
			&a'_{41}x_1 + a'_{42} x_2
		\end{bmatrix}^\top .
	\end{equation*}
	Then based on the definition of $\Delta^{ESS}$, we have
	\begin{equation*}
		x \in \Delta^{NE}
		\iff
		a'_{31}x_1 + a'_{32}x_2 \leq 0 \text{ and } a'_{41}x_1 + a'_{42} x_2 \leq 0.
	\end{equation*}
	Moreover, $a'_{41},a'_{31}>0$ in view of Lemma \ref{lem-signStructure}. 	
	So 
	\begin{align*}
		x \in \Delta^{NE}
		&\iff	
		x_1 + \frac{a'_{32}}{a'_{31}} x_2 \leq 0 \text{ and } x_1 + \frac{a'_{42}}{a'_{41}} x_2 \leq 0	\\
		&\iff	
		0\leq x_1 \leq \min\left\{ -\frac{a'_{32}}{a'_{31}},-\frac{a'_{42}}{a'_{41}} \right\}, x_1\leq 1.
	\end{align*}	
	Hence, if $\min\left\{ -\frac{a'_{32}}{a'_{31}},-\frac{a'_{42}}{a'_{41}} \right\} < 0$, then $x\not\in\Delta^{NE}$.
	Otherwise, $x = \alpha p^1+ (1-\alpha)p^2$ where 
	$
		\alpha \in \left[0, \min\left\{-\frac{a'_{32}}{a'_{31}}, - \frac{a'_{42}}{a'_{41}} ,1\right\} \right]
	$.
	Substituting the values of $a'_{ij}$ from $A'$ in the above equation completes the proof.
\end{IEEEproof}	

The following lemma reveals those singleton boundary equilibria that belong to $\Delta^{NE}$.
\begin{lemma}		\label{lem-singletonEquilibriaNash}
	$x^{13},x^{24}\not\in\Delta^{NE}$ and $x^{14}\in\Delta^{NE}$.
	Moreover, if $\S<\R<\frac{\T+\S}{2}$, or $m=2n+1,n\geq 1$ and $\R= \frac{\T+\S}{2}$, then $x^{23}\in\Delta^{NE}$.
	Otherwise, $x^{23}\not\in\Delta^{NE}$.
\end{lemma}
\begin{IEEEproof}
	The sign-structure of $A'x^{13}$ is of the form $\begin{bmatrix} + &+ &+ &++\end{bmatrix}^\top $. 
	Hence, $(p^4)^\top A'x^{13}>{x^{13}}^\top A'x^{13}$. 
	Hence, $x^{13}\not\in\Delta^{NE}$ by definition.
	Similarly $x^{24}\not\in\Delta^{NE}$ can be shown. 
	Now the result for $x^{23}$ is proven and that for $x^{14}$ can be done similarly.
	Define
	$
		z :=
		A'x^{23} = 
		\begin{bmatrix}
			 a'_{13}x_3 	
			&a'_{23}x_3	
			&a'_{32}x_2 
			&a'_{42}x_2
		\end{bmatrix}^\top .
	$
	Let $\S<\R<\frac{\T+\S}{2}$ or $m=2n+1,n\geq1$ and $\R=\frac{\T+\S}{2}$.
	In view of Lemma \ref{lem-signStructure}, $a'_{32}>a'_{42}$ and hence $z_3>z_4$.
	Similarly, $z_2\geq z_1$.
	Moreover, it can be verified that $z_2=z_3$.
	Hence, $z_2,z_3 = \max_{i\in\{1,\ldots,4\}} z_i$.
	Hence, any $x\in \Delta(p^2,p^3)$, maximizes $x^\top  z = x^\top A'x^{23}$ over $\Delta$. 
	Hence, since $x^{23}\in \Delta(p^2,p^3)$, it holds that ${x^{23}}^\top A'x^{23} \geq y^\top A'x^{23}$ for all $y\in\Delta$.
	Hence, $x^{23}\in\Delta^{NE}$.
	For all other payoffs, either $x^{23}\not\in\Delta$ or $z_1>z_2$.
	The first case clearly implies $x^{23}\not\in\Delta^{NE}$.
	For the second case, $(p^1)^\top A'x^{23}>{x^{23}}^\top A'x^{23} $, which rules out $x^{23}$ from $\Delta^{NE}$.
\end{IEEEproof}

%%%
\subsection{Convergence to a line segment implies convergence to a (set of continuum) stationary point(s)} \label{app:line}
In the analysis of Section \ref{subsec-interiorTrajectories}, we often face the situation where we know that the trajectory converges to a line segment. 
However, for completeness of our convergence results, we need to know whether the omega limit set of the trajectory is the whole line segment or just some parts of it.
For this purpose, we use a theorem showing that if a trajectory converges to a line segment, it converges to an equilibrium point or a continuum of equilibria on that line segment. 
Consider the function $y:\mathbb{R}\to \mathbb{R}^n$ and the set $\str S\subseteq \mathbb{R}^n$.
The expression $y(t)\to\str S$ implies that for any $\epsilon>0$, there exists some $M>0$ such that 
$
	t>M \Rightarrow \inf_{s\in \str S} \| y(t) - s\| < \epsilon
$
where $\|\cdot\|$ denotes an arbitrary norm in $\mathbb{R}^n$.
\begin{theorem}[reformulation of Corollary 1 in \cite{ramazi2017limit1}]			\label{lem-wLimitOfLine}
	Consider the $\mathbf{C}^r,r\geq 1,$ vector field 
	%\begin{equation}	\label{lem-wLimitOfLine-0}
	$	\dot{y} = f(y) , y\in\mathbb{R}^n.$
	%\end{equation}	
	If $y(t)$ converges to a compact simple open curve $\str L$, then $y(t)$ converges to an equilibrium point or a continuum of equilibrium points on $\str L$.
\end{theorem}

\section*{Acknowledgments}
	We would like to thank Dr. Hildeberto Jard\'on-Kojakhmetov for his technical discussions.

   \bibliographystyle{IEEEtran}
   \bibliography{../bib.bib}
%%%%%%%%%%%%%%%%%%%%%%%%%%%%%%%%%%%%%%%%
\begin{IEEEbiography}[{\includegraphics[width=1in,height=1.25in,clip,keepaspectratio]{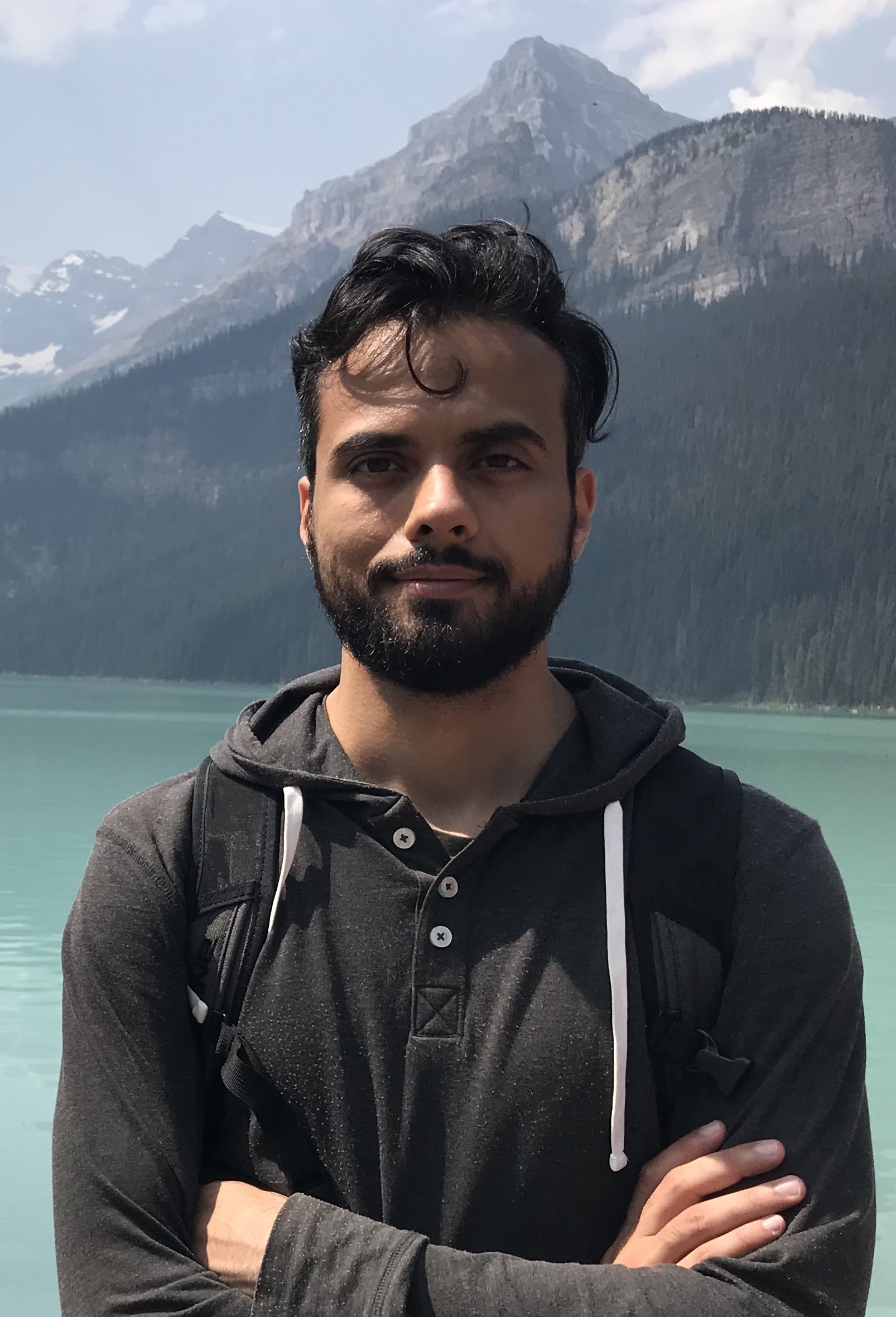}}]{Pouria Ramazi}
	 received the B.S. degree in electrical engineering in 2010 from University of Tehran, Iran, the M.S. degree in systems, control and robotics in 2012 from Royal Institute of Technology, Sweden, and the Ph.D. degree in systems and control in 2017 from the University of Groningen, the Netherlands. 
He is currently a joint Postdoctoral Research Associate with the Departments of Mathematical and Statistical Sciences and Computing Sciences of the University of Alberta.
\end{IEEEbiography}

\begin{IEEEbiography}[{\includegraphics[width=1in,height=1.25in,clip,keepaspectratio]{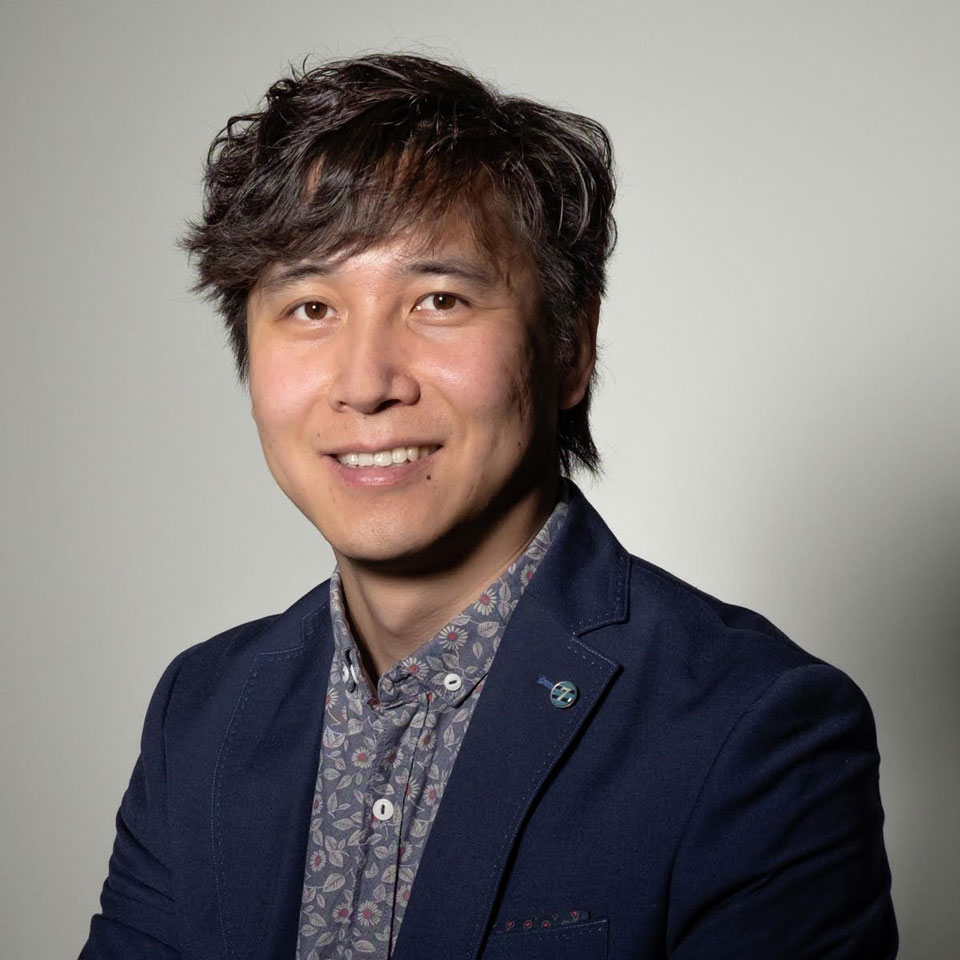}}]{Ming Cao}
	 is currently a professor of systems and control with the Engineering and Technology Institute (ENTEG) at the University of Groningen, the Netherlands, where he started as a tenure-track assistant professor in 2008. He received the Bachelor degree in 1999 and the Master degree in 2002 from Tsinghua University, Beijing, China, and the PhD degree in 2007 from Yale University, New Haven, CT, USA, all in electrical engineering. From September 2007 to August 2008, he was a postdoctoral research associate with the Department of Mechanical and Aerospace Engineering at Princeton University, Princeton, NJ, USA. He worked as a research intern during the summer of 2006 with the Mathematical Sciences Department at the IBM T. J. Watson Research Center, NY, USA. He is the 2017 and inaugural recipient of the Manfred Thoma medal from the International Federation of Automatic Control (IFAC) and the 2016 recipient of the European Control Award sponsored by the European Control Association (EUCA).  He is an associate editor for IEEE Transactions on Automatic Control, IEEE Transactions on Circuits and Systems and Systems \& Control Letters, and for the Conference Editorial Board of the IEEE Control Systems Society. He is also a member of the IFAC Technical Committee on Networked Systems. His main research interest is in autonomous agents and multi-agent systems, mobile sensor networks and complex networks. 
\end{IEEEbiography}   
   
   \end{document}